\documentclass{birkjour}
\usepackage{amsmath}
\usepackage{amssymb}

\usepackage{url}

\usepackage{yhmath}
\usepackage[dvips]{color}

\usepackage{srctex} 

\def\m#1{\mathsf{#1}} 

\newcommand{\cl}[2]{\ensuremath{\mathit{Cl}_{#1,#2}}}

\DeclareMathOperator{\Det}{Det} 

\DeclareMathOperator{\arccot}{arccot}

\DeclareMathOperator{\arccosh}{arcosh}
\DeclareMathOperator{\arcsinh}{arsinh}
\DeclareMathOperator{\arctanh}{artanh}

\DeclareMathOperator{\arccoth}{arcoth}

\newcommand{\bbZ}{\ensuremath{\mathbb{Z}}}
\newcommand{\bbN}{\ensuremath{\mathbb{N}}}

\newcommand{\reverse}[1]{\widetilde{#1}}
\newcommand{\gradeinverse}[1]{\wideparen{#1}}
\newcommand{\cliffordconjugate}[1]{\widetilde{\wideparen{#1}}}

\newcommand{\ii}{\mathrm{i}}

\newcommand{\ee}{\mathrm{e}} 

\def\A{\mathsf{A}}

\def\m#1{\mathsf{#1}}

\def\e#1{\mathbf{e}_{#1}} 

\newcommand{\ba}{\ensuremath{\mathbf{a}}}
\newcommand{\bb}{\ensuremath{\mathbf{b}}}

\newcommand{\bu}{\ensuremath{\mathbf{u}}}

\newcommand{\cA}{\ensuremath{\mathcal{A}}}
\newcommand{\cB}{\ensuremath{\mathcal{B}}}

\newcommand{\cU}{\ensuremath{\mathcal{U}}}
\newcommand{\cF}{\ensuremath{\mathcal{F}}}

\renewcommand{\d}{\ensuremath{\mathbin{\cdot}}} 
\newcommand{\w}{\ensuremath{\mathbin{\wedge}}} 

\usepackage{color}

\newcommand{\mycomment}[1]{} 

\newtheorem{thm}{Theorem}[section]

 \theoremstyle{definition}
 
 \theoremstyle{remark}

 \numberwithin{equation}{section}

\newcommand{\magnitude}[1]{\lvert #1\rvert}

\begin{document}


\title[Logarithm of multivector in 3D geometric algebras]
{Logarithm of multivector in real 3D Clifford algebras}

\author{A.~Acus}
\address{Institute of Theoretical Physics and Astronomy,\br
Vilnius University,\br Saul{\.e}tekio 3, LT-10257 Vilnius,
Lithuania} \email{arturas.acus@tfai.vu.lt}
\thanks{$\dagger$ Corresponding author: A.~Acus}

\author{A.~Dargys}
\address{%
Center for Physical Sciences and Technology,
\br Saul{\.e}tekio 3, LT-10257 Vilnius,
Lithuania} \email{adolfas.dargys@ftmc.lt}

\subjclass{Primary 15A18; Secondary 15A66}

\keywords{ {C}lifford (geometric) algebra, logarithms of Clifford
numbers, computer-aided theory}


\begin{abstract}

Closed form expressions for a logarithm  of general multivector
(MV) in base-free form  in real geometric algebras (GAs)
$\cl{p}{q}$ are presented for all $n=p+q=3$. In contrast to
logarithm of complex numbers (isomorphic to \cl{0}{1}), 3D
logarithmic functions, due to appearance of two double angle arc
tangent functions, allow to include \textit{two sets of sheets}
characterized by discrete coefficients. Formulas for generic and
special cases of individual blades and their combinations are
provided.
\end{abstract}

\maketitle





\section{\label{sec:intro}Introduction}

Logarithm properties are well-known for real and complex numbers.
Except the Hamilton quaternions which are isomorphic to \cl{0}{2},
the properties of logarithm in other 2D algebras (some partial
formulas for 2D GAs are provided in
\cite{Chappell2015,Josipovic2015,JosipovicErratum2015,Dargys2022Log})
and higher dimensional Clifford algebras remain uninvestigated as
yet. In general, GA logarithm properties are simplest for
anti-Euclidean  algebras
\cl{0}{n}~\cite{Gurlebeck1997,Pzreworska-Rolewicz1998}. As in the
complex algebra case we expect at least to have a principal
logarithm and a part that makes the GA logarithm a multivalued
function.

Recently in papers~\cite{Dargys2022a,Acus2022}, which will be the
starting point for the present article, we have performed a
detailed investigation of 3D exponential functions in real GAs.
However, the GA logarithm is more difficult to analyze since one
must take into account the multi-valuedness and the fact that in
3D algebras (except \cl{0}{3}) the logarithm may not exist for all
MVs. Here, we have treated the logarithm  as an inverse problem
using for this purpose the \textit{Mathematica} symbolic package,
more precisely as an inverse GA function to exponential in
separate 3D algebras \cl{0}{3}, \cl{3}{0}, \cl{1}{2}, and
\cl{2}{1}. The final GA logarithm formulas were checked
symbolically as well as numerically. They are in complete
agreement with more general formulas~\cite{AcusLNCS2023} suitable
for computation on any function of diagonalizable multivector
(MV). The exact logarithm formulas also have been applied to study
convergence of series expansion MV logarithms.

In Sec.~\ref{sec:notations} the notation is introduced. Since the
logarithm is closely related with a two argument arc tangent
function $\arctan(x,y)$, its properties are summarized in this
section as well. In Sec.~\ref{sec:Cl03} the logarithm of the
simplest, namely \cl{0}{3} algebra is considered. Since algebras
\cl{3}{0} and \cl{1}{2} are isomorphic, in Sec.~\ref{sec:Cl30} the
logarithms for both algebras are investigated simultaneously. In
Sec.~\ref{sec:Cl21} the logarithm of the most difficult \cl{2}{1}
algebra is presented. Since GA logarithm may be applied to GA
square root
calculation~\cite{Dargys2022Log,AcusDargysPreprint2020}, in
Sec.~\ref{sec:ArbitraryPowers} the roots and arbitrary fractional
powers of MV are discussed. In Sec.~\ref{sec:relations} the
relations of the logarithm to GA inverse trigonometric and
hyperbolic functions are presented. Finally, in
Sec.~\ref{sec:discussion} we summarize the obtained results.

\section{Notation and general properties of GA  logarithm}
\label{sec:notations} A general MV in GA space  is expanded  in
the orthonormal basis in inverse degree lexicographic ordering:
$\{1,\e{1},\e{2},\e{3},\e{12},\e{13},\e{23},\e{123}\equiv I\}$,
where $\e{i}$ are basis vectors, $\e{ij}$ are the bivectors and
$I$ is the pseudoscalar. The number of subscripts indicates the
grade. The scalar is a grade-0 element, the vectors $\e{i}$ are
the grade-1 elements, etc. In the orthonormalized basis the
geometric products of basis vectors satisfy the anticommutation
relation, $\e{i}\e{j}+\e{j}\e{i}=\pm 2\delta_{ij}$. For \cl{3}{0}
and \cl{0}{3} algebras the squares of basis vectors,
correspondingly, are $\e{i}^2=+1$ and $\e{i}^2=-1$, where
$i=1,2,3$. For mixed signature algebras such as  \cl{2}{1} and
\cl{1}{2} the squares are $\e{1}^2=\e{2}^2=1$, $\e{3}^2=-1$ and
$\e{1}^2=1$, $\e{2}^2=\e{3}^2=-1$, respectively.

The general MV in real Clifford algebras $\cl{p}{q}$ for $n=p+q=
3$ is
\begin{equation}\begin{split}\label{mvA}
\A=&\,a_0+a_1\e{1}+a_2\e{2}+a_3\e{3}+a_{12}\e{12}+a_{23}\e{23}+a_{13}\e{13}+a_{123}I\\
  \equiv&\,a_0+\ba+\cA+a_{123}I = a_0+\m{A}_{1,2}+a_{123}I=\m{A}_{0,1,2,3},
\end{split}
\end{equation}
where $a_i$, $a_{ij}$ and $a_{123}$ are the real coefficients, and
$\ba=a_1\e{1}+a_2\e{2}+a_3\e{3}$ and
$\cA=a_{12}\e{12}+a_{23}\e{23}+a_{13}\e{13}$ is, respectively, the
vector and bivector. $I$~is the pseudoscalar, $I=\e{123}$. The
double index in MV $\m{A_{i,j}}$ indicates the sum of MVs of grades
$i$ and $j$, i.e.
$\m{A_{i,j}}=\langle\m{A}\rangle_{i}+\langle\m{A}\rangle_{j}$.

The main involutions, namely the reversion, grade inversion and
Clifford conjugation denoted, respectively, by tilde, circumflex
and their combination are defined by
\begin{equation}\begin{split}
&\widetilde{\m{A}}=a_0+\ba-\cA-a_{123}I,\quad \gradeinverse{\m{A}}=a_0-\ba+\cA-a_{123}I,\\
&\cliffordconjugate{\m{A}}=a_0-\ba-\cA+a_{123}I.
\end{split}
\end{equation}

\subsection{\label{expLogTan} General properties of GA logarithm}

The logarithm of MV is another MV that  belongs to the same
geometric algebra (GA). The defining equation for MV logarithm is
$\log(\ee^{\m{A}})=\m{A}$, where $\m{A}\in\cl{p}{q}$. The GA
logarithm is a multivalued function. In \cite{Chappell2015} it was
suggested  that "The principal value of the logarithm can be
defined as the MV $\m{M}=\log (\m{Y})$ with the smallest norm",
where $Y\in\cl{p}{q}$. The natural norm for a MV is the
determinant norm defined in subsection~\ref{logSeriesSection}. The
following properties hold for MV logarithm:
\begin{equation}\label{logFormulas}\begin{split}
&\log(\m{A}\m{B})=\log(\m{A})+\log(\m{B})\quad\text{if\ } \m{A}\m{B}=\m{B}\m{A},\\
&\ee^{\log(\m{A})}=\m{A},\quad \ee^{-\log(\m{A})}=\m{A}^{-1}, \\
&\widetilde{\log(\m{A})}=\log(\widetilde{\m{A}}),
\quad\gradeinverse{\log(\m{A})}=\log(\gradeinverse{\m{A}}),
\quad
\cliffordconjugate{\log(\m{A})}=\log(\cliffordconjugate{\m{A}}),\\
&\m{V}\,\log(\m{A})\m{V}^{-1}=\log(\m{V}\m{A}\m{V}^{-1}).
\end{split}\end{equation}
In the last expression the transformation $\m{V}$, for example the
rotor, is  pushed inside the logarithm.

\subsection{\label{logSeriesSection} GA logarithm series}
In analogy with a definition of logarithm in  complex plane for GA
logarithm  we can write
\begin{equation}\label{logSeries}
  \log \m{A} = \sum_{k=1}^{\infty} \frac{(-1)^{k - 1} (\m{A} - 1)^k}{k}, \quad \mathrm{if}\quad |\m{A} - 1| < 1,
\end{equation}
where $|\m{A} - 1|$ denotes the determinant norm.  For arbitrary
MV the determinant norm is defined as an absolute value of
determinant $\Det(\m{B})$ of MV~$\m{B}$ raised to fractional power
$1/k$, where $k=2^{\lceil n/2\rceil}$, i.e.,
$|\m{B}|=\bigr(\Det(\m{B})\bigl)^{1/k}>0$. For algebras having
negative determinant instead the semi-norm (aka pseudoscalar) is
introduced
$\magnitude{\m{B}}=\bigr(\text{abs}(\Det(\m{B}))\bigl)^{1/k}\ge
0$. The equality sign means that in case of semi-norm the
determinant may be zero although $\m{B}\ne 0$.  In the
following the same symbol will be used for both the norm
and semi-norm. The (semi-)norm can be interpreted as a number of
multipliers needed to define $\Det(\m{B})$. In 3D algebras ($n=3$)
we have $k=2^{\lceil 3/2\rceil}=2^{2}=4$, which is the degree of
characteristic~\cite{Abdulkhaev2021} polynomial $\Det(\m{B})$.
In this way found integer  $k$ coincides with the number of
multipliers in the 3D determinant: $\Det(\m{B})= \m{B}
\reverse{\m{B}} \gradeinverse{\m{B}}
\gradeinverse{\reverse{\m{B}}}$. The determinant norm for MV
$\m{B}$ in 3D algebras, therefore, is
$|\m{B}|=\sqrt[4]{\mathop{\mathrm{abs}}(\Det(\m{B}))}$\,. It can
be shown that for any GA that holds a basis element with
property $\e{i}^2=-1$ by adding a scalar one can construct a MV
the norm of which may be identified with a module of a complex
number. For example in \cl{3}{0} the norm of $\m{B}=1+\e{12}$ is
$\sqrt{(1+\e{12})(1-\e{12})}=\sqrt{2}$ which coincides with
$|\m{B}|=\sqrt[4]{\mathop{\mathrm{abs}}(\Det(\m{B}))}=\sqrt{2}$\,.
(also refer to Example~1 below).

If the MV has a numerical form, to minimize the number of
multiplications it is convenient to represent the logarithm  in a
nested form (aka Horner's rule). The logarithmic
series~\cite{Abramowitz1964} (also called Mercator series), if
rewritten according to Horner's rule,  assumes the following form,
\begin{equation}\begin{split}\label{logHorner}
\log\m{B}=\m{B}(1+\m{B}(-\tfrac12+\m{B}(\tfrac13+\m{B}
(-\tfrac14+\m{B}(\tfrac15+\dotsm ))))),\quad\textrm{where\ }
\m{B}=\m{A}-1.
\end{split}\end{equation}
\vspace{3mm}
 \textbf{Example~1.}
{\it  MV equivalent to complex number.}  Let's take the MV
$\m{A}=\frac{9}{10}-\frac{1}{3}\e{3}$  the determinant norm of
which in \cl{0}{3} is
$|\m{B}|=|\m{A}-1|=\frac{\sqrt{109}}{30}\approx0.34801<1$.
Therefore, the standard series, Eq.~\eqref{logSeries}, may be
applied to find an approximate value (the result found by exact
formula in Example~2 is $\log \frac{\sqrt{829}}{30}-
\e{3}\arctan\frac{10}{27}\approx-0.0410873 - \e{3}0.354706$).
Since $\e{3}^2=-1$ and it is the only basis vector in the
considered MV, one may replace the MV by complex number
$z=\frac{9}{10}-\ii\frac{1}{3}$. The module is
$|z-1|=\frac{\sqrt{109}}{30}$ which  coincides with the MV
determinant norm. Then, $\log z \approx-0.0410873 - \ii\,
0.354706$.

Now let's calculate the logarithm of
$\m{A}^\prime=-\frac{9}{10}-\frac{1}{3}\e{3}$ by the Horner
series~\eqref{logHorner}. Since
$|\m{A}^{\prime}-1|=\frac{\sqrt{3349}}{30} \approx 1.92902>1$ the
series diverges. As shown in Example~2, the logarithm can be
easily computed if exact GA logarithm formula obtained in
the present paper is used. After replacement of the MV by complex
number we obtain that $|z^{\prime}-1|=1.92902$ which again
coincides with the module of
$z^{\prime}=-\frac{9}{10}-\frac{1}{3}\ii$. Computing the value of
logarithm by \textit{Mathematica} command
$\mathrm{FunctionExpand}[\mathrm{Log}[z^{\prime}]]$ we obtain
$\log z^{\prime}=\ii (-\pi+\arctan(10/27))-\log(30)+\log(829)/2$
which has the same numerical value as shown in  the Example~2.

\subsection{Double-argument arc tangent function}\label{arctanProperties}
\begin{figure}[t]
\centering
\includegraphics[width=7cm]{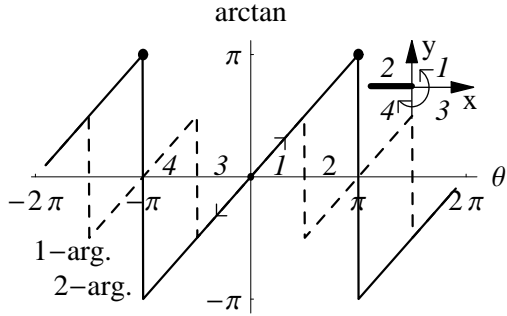}
\caption{Graphical representation of single
$\arctan(y/x)=\arctan(\sin\theta/\cos\theta)$ (dashed line) and
double  $\arctan(x,y)=\arctan(\cos\theta,\sin\theta)$ (solid line)
argument tangent functions used by \textit{Mathematica}. $\theta$
is an angle between $x$ axis and vector (not shown) attached to
the center of complex $x-y$ plane. The vector may be rotated from
$x$-axis anticlockwise, $\theta=(0...\pi]$, or clockwise,
$\theta=[0...-\pi)$ as shown by arrows  in the inset. In the inset
also the numbering of the quadrants $\textit{1-4}$ and  the
branching, represented by thick line on the negative part of $x$
axis, are shown.\label{fig:arctan}}
\end{figure}

GA logarithm as we shall see, in its nature  is a multi-valued
function with period $2\pi$. To account for quadrant sign in
complex plane properly we shall need the double argument arc
tangent function as given in the \textit{Mathematica}, the
properties of which are briefly mentioned below.
Figure~\ref{fig:arctan} shows the  single  and double argument arc
tangent functions. The former has period $\pi$ and its principal
values lie in the interval $\theta=[-\pi/2,\pi/2)$, while the
double argument arc tangent has, respectively, $2\pi$ period  and
principal values in $\theta=[-\pi,\pi)$. The inset on the right
side  of Fig.~\ref{fig:arctan} shows the quadrants $\textit{1-4}$
in  $x-y$ plane. Note that the anticlockwise rotation is done from
quadrant \textit{1} to quadrant \textit{2}, while clockwise
rotation in order \textit{3$\to$4}, so that a jump in the double
arc tangent value and associated branching occurs on the negative
side of $x$-axis rather than on $y$ axis as is in the standard
single argument case. Also, in Fig.~\ref{fig:arctan} note small
points on vertical branching steps that indicate that respective
arc tangent value on periodic line belong to upper rather than
lower part, i.e., at $\theta=\pi$ we have
$\arctan(\cos(\pi),\sin(\pi))=\arctan(-1,0)=\pi$, however, after
addition of infinitesimal  angle
$\arctan(\cos(\pi+0_{+}),\sin(\pi+0_{+}))=-\pi$. Similarly, at
$\theta=-\pi$ we have $\arctan(\cos(-\pi),\sin(-\pi))=\pi$, and,
$\arctan(\cos(-\pi+0_{+}),\sin(-\pi+0_{+}))=-\pi$. If $x,y$ were
replaced by real numbers \textit{Mathematica} will switch
automatically to a single argument arc tangent in the first
quadrant and principal values, for example,
$\arctan(17,10)=\arctan(10/17)$,
$\arctan(-17,10)=\pi-\arctan(10/17)$,
$\arctan(17,-10)=-\arctan(10/17)$,
$\arctan(-17,-10)=-\pi+\arctan(10/17)$.

In the terms of a standard arc tangent function the argument of
which is in the range $(-\pi/2, \pi/2)$, the double tangent
principal values now in the range $(-\pi, \pi)$ can be expressed
as follows:
\begin{equation}
  \arctan(x,y)=\begin{cases}\arctan({\frac {y}{x}})&{\text{if }}x>0,\\
    \arctan({\frac {y}{x}})-\pi &{\text{if }}x<0{\text{ and }}y\geq 0,\\
    \arctan({\frac {y}{x}})+\pi &{\text{if }}x<0{\text{ and }}y<0,\\
    +{\frac {\pi }{2}}&{\text{if }}x=0{\text{ and }}y>0,\\
    -{\frac {\pi }{2}}&{\text{if }}x=0{\text{ and }}y<0,\\
    \textrm{undefined}&{\text{if }}x=0{\text{ and }}y=0\, .
  \end{cases}
\end{equation}

We will start from $\cl{0}{3}$  where the expanded exponential in
a coordinate-form has the simplest MV coefficients and the
logarithm exists for all MVs.

\section{\label{sec:Cl03} MV  logarithms in \cl{0}{3}}
\vspace{3mm}

\subsection{\label{allCasesCl03} Logarithm formula for generic MV}
The term "generic" here will be understood as "not creating the
problems". If for a given set of MV coefficients the generic
formula is not applicable, for example, due to nullification of a
denominator, or due to appearance of an undefined subexpression
like $\arctan(0,0)$, we will refer to it as "special case".
Special cases will be covered by more elaborate formula later.

\begin{thm}[Logarithm of multivector in
\cl{0}{3}]\label{theoremCl03}  The generic logarithm of MV
$\m{A}=a_0+\ba+\cA+a_{123}I$ is the MV given by
\begin{equation}\label{log03coordfree}
\log(\m{A})= \frac{1}{2}\bigl(\m{A}_{{0}_{+}}+\m{A}_{{0}_{-}}+\m{A}_{{1,2}_{+}}+\m{A}_{{1,2}_{-}}+(\m{A}_{{0}_{+}} -\m{A}_{{0}_{-}})I\bigr),
\end{equation}
  with
  \begin{align}
    \m{A}_{{0}_{+}}&=\frac{1}{2}\log\bigl((a_{0}+a_{123})^2+a_{+}^2\bigr),&a_{+}\neq 0,\label{generic03allterms1}\\
    \m{A}_{{0}_{-}}&=\frac{1}{2}\log\bigl((a_{0}-a_{123})^2+a_{-}^2\bigr),&a_{-}\neq 0,\label{generic03allterms2}\\
  \m{A}_{{1,2}_{+}}&=\frac{1}{a_{+}}\big(\arctan (a_{0}+a_{123},a_{+})+
    2\pi c_{1+}\big)\bigl(1+I\bigr)\bigl(\ba +\cA\bigr),&a_{+}\neq 0,\label{generic03allterms3}\\
  \m{A}_{{1,2}_{-}}&=\frac{1}{a_{-}}\big(\arctan
(a_{0}-a_{123},a_{-})+2\pi c_{1-}\big)\bigl(1-I\bigr)\bigl(\ba
    +\cA\bigr),&a_{-}\neq 0.\label{generic03allterms4}
\end{align}
  The MVs $\m{A}_{{0}_{\pm}},\m{A}_{{1,2}_{\pm}}$ and $\m{A}_{{0}_{\pm}} I$ denote,
respectively,  the scalar, vector$\pm$bivector and the
pseudoscalar components.  $c_{1_{\pm}},c_{2_{\pm}}\in \mathbb{Z}$
  are arbitrary integers. The scalars $a_{+}\ge 0$ and $a_{-}\ge 0$
are given by expressions~\cite{Dargys2022a,Acus2022},
\begin{align}
  a_{-}=&\sqrt{-(\ba\d\ba+\cA\d\cA) + 2I  \ba\w\cA}\label{aminusCl03}\\
  =&\sqrt{(a_{3}+a_{12})^2+(a_{2}-a_{13})^2+(a_{1}+a_{23})^2},\notag\\
  a_{+}=&\sqrt{-(\ba\d\ba+\cA\d\cA) - 2I  \ba\w\cA}\label{apliusCl03}\\
  =&\sqrt{(a_{3}-a_{12})^2+(a_{2}+a_{13})^2+(a_{1}-a_{23})^2}\,,\notag
\end{align}
\end{thm}

\textit{Proof}. It is enough to check that after substitution of
\eqref{log03coordfree} into GA exponential formula (1)~of~\cite{Dargys2022a} one will get the initial MV~$\m{A}$.

The Theorem~\ref{theoremCl03} gives the GA logarithm  in a
basis-free form. However, the derivation of above given
generic logarithm formula at first was done in a coordinate
form, from which the Theorem~\ref{theoremCl03} follows (see
Appendix~\ref{coordFormCL03}). The Theorem~\ref{theoremCl03}
ensures the existence of GA logarithm for all MVs with real
coefficients in $\cl{0}{3}$, because in the mentioned  algebra the
zero determinant of MV ($\Det\m{A}=0$) occurs only if $\m{A}=0$.
As we shall see this property does not hold for remaining
algebras.

\subsection{\label{spacialCL04} Special cases}

In Theorem~\ref{theoremCl03} it was presumed that  the both
scalars $a_{-}$ and $a_{+}$ do not vanish. This assumption is
equivalent to the condition that the sum of vector and bivector
must have non-zero-determinant\footnote{We shall remind that for
$\cl{0}{3}$ algebra the determinant  $\Det\m{A}=0$ means
$\m{A}=0$. This property makes the exponential/logarithm
analysis in this algebra relatively simple. In remaining 3D
algebras additional conditions are needed to include MVs with
$\Det\m{A}=0$.}, $\Det(\ba+\cA)= a_{+}^2 a_{-}^2\neq 0$. If either
of scalars is zero then we have a special case. This situation is
met in rare cases, for instance\footnote{Minus sign in
$a_2=-a_{13}$ comes from a strictly increasing order of numbers in
the basis element indices, due to the so-called inverse degree
lexicographic ordering~\cite{Dargys2022a}.}, when $a_1=a_{23}$,
$a_2=-a_{13}$, $a_3=a_{12}$. In such and similar cases the MVs
$\m{A}_{{0}_{\pm}},\m{A}_{{1,2}_{\pm}}$ and $\m{A}_{{0}_{\pm}} I$
in the Theorem~\ref{theoremCl03} must be supplemented by conditions:
  \begin{align}
  \m{A}_{{0}_{+}}=&\begin{cases}
    \log\bigl(a_{0}+ a_{123}\bigr)+  2\pi c_{2_{+}}\hat{\cU},& a_{+}=0\quad \textrm{and}\quad (a_{0}+ a_{123}) > 0\\
    \log(0_{+}),&   a_{+}=0\quad \textrm{and}\quad (a_{0}+ a_{123}) = 0\\
    \log\bigl(-(a_{0}+ a_{123})\bigr)\\
\qquad
    +  (\pi+2\pi c_{2_{+}})\hat{\bu},& a_{+}=0\quad \textrm{and}\quad (a_{0}+ a_{123}) < 0,\label{B0ForCl03FullAplus}
\end{cases}\\[10pt]
\m{A}_{{0}_{-}}=&\begin{cases}
    \log\bigl(a_{0}- a_{123})+  2\pi c_{2_{-}}\hat{\cU},& a_{-}=0\quad \textrm{and}\quad (a_{0}- a_{123}) > 0\\
    \log(0_{+}),&   a_{-}=0\quad \textrm{and}\quad (a_{0}- a_{123}) = 0 \\
    \log\bigl(-(a_{0}- a_{123})\bigr)\\
\qquad
    +  (\pi+2\pi c_{2_{-}})\hat{\bu},& a_{-}=0\quad \textrm{and}\quad (a_{0}-a_{123}) < 0,\label{B0ForCl03FullAminus}
\end{cases}\\
%
  \m{A}_{{1,2}_{+}}=&\begin{cases}
    \bigr(\frac{1}{a_{0}+ a_{123}} +2\pi c_{1_{+}} \bigr)\\
\qquad \times  (1+I)\bigl(\ba +\cA\bigr),& a_{+}=0\quad \textrm{and}\quad (a_{0}+ a_{123}) > 0\\
0,& a_{+}=0\quad \textrm{and}\quad (a_{0}+a_{123})= 0,\\
    (\pi+2\pi c_{1_{+}}) (1+I)\bigl(\ba +\cA\bigr),& a_{+}=0\quad \textrm{and}\quad (a_{0}+a_{123}) < 0,\label{B12ForCl03FullPlus} \\
\end{cases}\\[10pt]
  \m{A}_{{1,2}_{-}}=&\begin{cases}
    \bigr(\frac{1}{a_{0}- a_{123}} +2\pi c_{1_{-}} \bigr)\\
\qquad \times  (1- I)\bigl(\ba +\cA\bigr),& a_{-}=0\quad \textrm{and}\quad (a_{0}-a_{123}) > 0\\
0,& a_{-}=0\quad \textrm{and}\quad (a_{0}- a_{123})= 0 \\
    (\pi+2\pi c_{1_{-}}) (1-I)\bigl(\ba +\cA\bigr),& a_{-}=0\quad \textrm{and}\quad (a_{0}-a_{123}) < 0,\label{B12ForCl03FullMinus}
\end{cases}
\end{align}
Here $c_{1_{\pm}},c_{2_{\pm}}\in \mathbb{Z}$ are the arbitrary
integers. The conditions  for $(a_0\pm a_{123})$ on the
right-hand side take into account the case $\Det(\ba+\cA)=0$. In
scalars\footnote{The appearance of a free vector/bivector
breaks the grade arrangement in the generic terms
\eqref{generic03allterms1} and \eqref{generic03allterms2}. The
choice, however, results in a more simple final expression, since
now it is enough to write only a single free $\bu$
or $\cU$ term (see Eq.~\eqref{log03coordfree}) instead of a pair
$\bu$ and $\bu I$, or $\cU$ and $\cU I$, if we would have chosen
to move these terms to vector+bivector part by following up a
strict grade notation convention.} $\m{A}_{{0}_{+}}$ and
$\m{A}_{{0}_{-}}$, the symbols  $\hat{\bu}$ and $\hat{\cU}$
represent any free unit vector or bivector, respectively,
$\hat{\bu}^2=\hat{\cU}^2=-1$. For example, the unit vector can be
parametrized as
$\hat{\bu}=(u_1\e{1}+u_2\e{2}+u_3\e{3})/\sqrt{u_1^2+u_2^2+u_2^3}$.
It should be noted that the term $1/(a_{0}\pm a_{123})$ in
Eqs~\eqref{B12ForCl03FullPlus} and \eqref{B12ForCl03FullMinus}
represents the limit $\lim_{a_{\pm}\to 0} \arctan (a_{0}\pm
a_{123},a_{\pm})/{a_{\pm}}=1/(a_{0}\pm a_{123})$ which is valid
only when $a_{0}\pm a_{123}>0$. The notation of $\log(0_{+})$ in
expressions for $\m{A}_{{0}_{+}}$ and $\m{A}_{{0}_{-}}$\ is
explained in Example~6.

Interpretation of special conditions
\eqref{B0ForCl03FullAplus}-\eqref{B12ForCl03FullMinus} in terms of
the MV determinant~\cite{Acus2018,Helmstetter2019,Shirokov2020a}
becomes more evident if one remembers that the determinant of MV
$\m{A}$ in $\cl{0}{3}$ can be expressed in a form
$\Det(\m{A})=\bigl(a_{-}^2+(a_{0}-a_{123})^2\bigr)\bigl(a_{+}^2+(a_{0}+a_{123})^2\bigr)$,
whereas the condition $a_{\pm}=0$ is equivalent to
$\Det(\m{A}_{{12}_{\pm}})=\Det(\ba+\cA)= a_{+}^2 a_{-}^2$. All
special cases therefore occur if $\Det (\ba+\cA)=0$ and the
condition are described by $a_{0}\pm a_{123}\lesseqqgtr 0$. In
conclusion, the symbolic expression for logarithm, has three
special pieces (branches) $a_0\pm a_{123}\lesseqqgtr 0$ provided
the condition $a_{\pm}= 0$ is satisfied and the generic piece is
characterized by $a_{\pm}\neq 0$.

\subsection{\label{freeParameters} Multivaluedness and free multivector}
To include multivaluedness in GA logarithm we introduce a free
multivector $\m{F}$ by the following  defining
equation~\cite{Dargys2022Log}
\begin{equation}
\ee^{\log(\m{A})+\m{F}}=\ee^{\log(\m{A})}\ee^{\m{F}}=\ee^{\log(\m{A})},
\end{equation}
which implies two conditions the MV $\m{F}$ must satisfy: the
commutator $[\log(\m{A}),\m{F}]=0$  and $\ee^{\m{F}}=1$. As we
shall see, for remaining $n=3$ algebras the free MV $\m{F}$ will
play a similar role. One can check that the expression
\begin{align}\label{periodicityCL03}
  \m{F}=& \frac{\pi c_{1_{+}}}{a_{+}}(1+I)\bigl(\ba +\cA\bigr)+\frac{\pi c_{1_{-}}}{a_{-}}(1-I)\bigl(\ba +\cA\bigr)
\end{align}
satisfies $\ee^{\m{F}}=1$, and that for a generic MV $\m{A}$,
Eq.~\ref{log03coordfree}, the free term~\eqref{periodicityCL03}
commutes with $\log(\m{A})$, i.e. $[\log(\m{A}),\m{F}]=0$. The
integer constants $c_{1_{+}},c_{1_{-}}\in \mathbb{Z}$ in
Eq.~\eqref{periodicityCL03} add two free (discrete) parameters
that may be used to shift  the coefficients of vector and bivector
in $\log(\m{A})$ by some multiple of~$\pi$. The sum $(\ba +\cA)$
in~\eqref{periodicityCL03} constitute vector+bivector
part\footnote{In 3D algebras, the scalar and pseudoscalar belong
to algebra center and as a result they commute with all elements.}
of the original MV $\m{A}$, therefore $(\ba +\cA)$ automatically
commutes with $\m{A}$. As a result only discrete free coefficients
are possible in the logarithm generic formula. In special cases
(see Eqs~\eqref{logaPlusbPlusrelationsam} and
\eqref{logaPlusbPlusrelationsap} in the
Appendix~\ref{coordFormCL03}) the free MV $\m{F}$ may also contain
arbitrary unit vector $\hat{\bu}$ and/or unit bivector $\hat\cU$.
In such cases one can include two additional continuous parameters
interpreted as directions of $\hat{\bu}$ or $\hat\cU$.

Since $\arctan(x,y)$ has been defined in the range $(-\pi,\pi]$
(usually called the principal value or the main branch,
Fig.~\ref{fig:arctan}), we can add to it any multiple of
$2\pi$. Therefore, the plus/minus instances of $\arctan (a_{0}\pm
a_{123},a_{\pm})$  (see Eqs  \eqref{logaPlusbPlusrelationsam} and
\eqref{logaPlusbPlusrelationsap})
 were
replaced by more general expressions $\arctan
(a_{0}+a_{123},a_{+})+2\pi c_{1_{+}}$ and $\arctan
(a_{0}-a_{123},a_{-})+2\pi c_{1_{-}}$ in
Eqs~\eqref{generic03allterms3} and \eqref{generic03allterms4},
respectively, which takes into account the multivaluedness of the
argument. This explains the rationale behind the construction of
the free MVs for GA logarithm.

In \cite{Higham08} the notion of principal logarithm (also called
the principal value of logarithm) in case of matrices was
introduced. In \cite{Chappell2015} it was suggested that the
"logarithm principal value in GA can be defined as the MV
$\m{M}=\log\m{Y}$ with the smallest norm".
Formulas~\eqref{generic03allterms1}-\eqref{generic03allterms4} and
\eqref{B0ForCl03FullAplus}-\eqref{B12ForCl03FullMinus} might
suggest that we could obtain the principal logarithm values after
equating discrete free constants $c_{1_{\pm}},c_{2_{\pm}}$ to
zero. Unfortunately, extensive numerical checks revealed that this
is not always the case.

\vspace{3mm}
 \textbf{Example 2.} {\it Logarithm of simple MV in
\cl{0}{3}.}  For MV $\m{A}=\tfrac{9}{10}-\tfrac{1}{3}\e{3}$ in the
Example~1,  Eqs~\eqref{aminusCl03} and \eqref{apliusCl03} give
$a_{+}=a_{-}=1/3$. The MVs in~\eqref{log03coordfree} then are
$\m{A}_{0+}=\m{A}_{0-}=-\log(10/9)+\pi\e{3}$,
$\m{A}_{1,2\pm}=\frac{\pi}{3}(-\e{3}\pm\e{12})$. The logarithm
calculated by exact formula~\eqref{log03coordfree} is
$\log(\tfrac{9}{10}-\tfrac{1}{3}\e{3})=-\tfrac12\log\frac{900}{829}-
\arctan\left(\frac{10}{27}\right)\e{3}\approx-0.0410873-0.354706\e{3}$
which  coincides with result  of Example~1. Now, let's calculate
GA logarithm of $\m{A}^\prime=-\tfrac{9}{10}-\tfrac{1}{3}\e{3}$
that diverges when the series~\eqref{logHorner} is used.  With
exact formulas~\eqref{generic03allterms1}--\eqref{apliusCl03} we
find: $a_{+}=a_{-}=\frac{1}{3}$,
$\m{A}_{0+}=\m{A}_{0-}=-\frac{1}{2}\log(900/829)$,
$\m{A}_{1,2+}=\big(\pi-\arctan(10/27)\big)(\e{12}-\e{3})$,
$\m{A}_{1,2-}=\big(\arctan(10/27)-\pi\big)(\e{12}+\e{3})$. Then
Eq.~\eqref{log03coordfree} gives $\log(\m{A}^\prime)=-\frac{1}{2}
\log(\frac{900}{829})+(\arctan(\frac{10}{27})-\pi) \e{3}\approx
-0.0410873-2.78689\e{3}$. Exponentiation of the obtained logarithm
gives  initial MV,
$\exp\big(\log(\m{A}^\prime)\big)=\m{A}^\prime$. The result also
can be checked by complex logarithm, because the initial MV
consist of scalar and basis vector $\e{3}^2=-1$ only.

\vspace{3mm}
 \textbf{Example~3.} {\it Logarithm of generic MV in \cl{0}{3}.}
Let's  compute the logarithm of $\m{A}=-8-6 \e{2}-9 \e{3}+5
\e{12}-5 \e{13}+6 \e{23}-4 \e{123}$. Then, $a_{-}^2=53$ and
$a_{+}^2=353$. The
Eqs~\eqref{generic03allterms1}-\eqref{generic03allterms4} give
$\m{A}_{{0}_{+}}=\tfrac12\log (497)$,
$\m{A}_{{0}_{-}}=\tfrac12\log (69)$,
$\m{A}_{{1,2}_{+}}=(353)^{-1/2}\bigl(\pi
-\arctan\bigl(\frac{\sqrt{353}}{12}\bigr)+2 \pi  c_{1_{+}}\bigr)
(1+ I)\bigl(-6 \e{2}-9 \e{3}+5 \e{12}-5 \e{13}+6 \e{23}\bigr)$ and
$\m{A}_{{1,2}_{-}}=(53)^{-1/2}\bigl(\pi
-\arctan\left(\frac{\sqrt{53}}{4}\right)+2 \pi c_{1_{-}}\bigr)(1-
I)(-6 \e{2}-9 \e{3}+5 \e{12}-5 \e{13}+6 \e{23})$, where the free
term $\m{F}$, Eq.~\eqref{periodicityCL03}, has been included via
$c_{1_{\pm}}$. The logarithm is the sum  of all above listed MVs:
$\log(\m{A})=
\frac{1}{2}\bigl(\m{A}_{{0}_{+}}+\m{A}_{{0}_{-}}+\m{A}_{{1,2}_{+}}+\m{A}_{{1,2}_{-}}+(\m{A}_{{0}_{+}}
-\m{A}_{{0}_{-}})I\bigr)$. Using the
exponential~\cite{Dargys2022a} one can check that the numerical
logarithm $\log(\m{A})$ indeed yields the initial MV $\m{A}$ for
arbitrary integer constants $c_{1_{\pm}}$.

\vspace{3mm}
 \textbf{Example~4.} {\it Logarithm of MV when $a_{+}=0$ and
$a_{0}+a_{123}>0$}. The MV that satisfies  these  conditions is
$\m{A}=1+(3 \e{1}-2 \e{2}+\e{3})+(\e{12}+2 \e{13}+3 \e{23})+7
\e{123}=1+\ba+\cA+7 \e{123}$. Equation
 \eqref{aminusCl03}  gives $a_{-}=\sqrt{56}=2\sqrt{14}$ and
$a_{0}-a_{123}=-6$. Eqs~\eqref{B0ForCl03FullAplus},
\eqref{B12ForCl03FullPlus} give $\m{A}_{{0}_{+}}=\log 8+ 2\pi
c_{2_{+}}\hat{\cU}$, $\m{A}_{{1,2}_{+}}=\bigl(\frac{1}{8}+2 \pi
c_{1_{+}}\bigr) (1+I)(\ba+\e{3}+\cA)=0$. Then from
\eqref{generic03allterms2} and \eqref{generic03allterms4} we have
$\m{A}_{{0}_{-}}=\tfrac{1}{2}\log 92$ and
$\m{A}_{{1,2}_{-}}=\frac{\pi
-\arctan\bigl(\frac{\sqrt{14}}{3}\bigr)+2 \pi
c_{1_{-}}}{2\sqrt{14}}(1- I)(\ba+\cA)$. Finally, from
Eq.~\eqref{log03coordfree} $\log(\m{A})=\frac{1}{28} \Bigl(7
\bigl(\log 5888-\log \frac{23}{16} \e{123}\bigr)+\sqrt{14}
\bigl((2 c_{1_{-}}+1)\pi -\arctan\frac{\sqrt{14}}{3}\bigr)
(\ba+\cA)\Bigr)+ (1+I)\pi c_{2_{+}}\hat{\cU} $. After
exponentiation of $\m{A}$ the constants $c_{1_{-}}$ and
$c_{2_{+}}$ and bivector $\hat{\cU}$ simplify out.

\vspace{3mm}
 \textbf{Example~5.} {\it Logarithm of MV when
$a_{-}=0$ and $a_{0}-a_{123}<0$.} These conditions are satisfied
by $\m{A}=1+(-3 \e{1}+2 \e{2}-\e{3})+(\e{12}+2 \e{13}+3 \e{23})+7
\e{123}=1+\ba+\cA+7 \e{123}$. We have $a_{+}^2=56$,
$a_0+a_{123}=9$ and $a_{0}-a_{123}=-6$. Then,
Eq.~\eqref{B0ForCl03FullAminus} gives
 $\m{A}_{{0}_{-}}=\log 6+
(\pi+  2\pi c_{2_{-}})\hat{\bu}$. The
Eqs~\eqref{generic03allterms1} and \eqref{generic03allterms3} give
$\m{A}_{{0}_{+}}=\frac12\log120$,
$\m{A}_{{1,2}_{+}}=\frac{1}{2\sqrt{14}}\bigl(\arctan\bigl(\tfrac{1}{2}\sqrt{7/2}\bigr)+2\pi
c_{1_{+}}\bigr) (1+I)(\ba+\cA)$, and
Eq.~\eqref{B12ForCl03FullMinus} $\m{A}_{{1,2}_{-}}=(\pi+2\pi
c_{1_{-}})(1-I)(\ba+\cA)=0$. Finally, $\log(\m{A})=
\frac{1}{2}\bigl(\m{A}_{{0}_{+}}+\m{A}_{{0}_{-}}+\m{A}_{{1,2}_{+}}+(\m{A}_{{0}_{+}}
-\m{A}_{{0}_{-}})I\bigr)$.

\vspace{3mm}
 \textbf{Example~6.} {\it Logarithm with infinite
subparts: The case $a_{+}=0$ and $a_{0}+a_{123}=0$.} The example
exhibits unusual and the most interesting case. In $\cl{0}{3}$,
let's  compute GA logarithm of $\m{A}=1+(-2 \e{1}-3 \e{2}+5
\e{3})+(5 \e{12}+3 \e{13}-2 \e{23})-\e{123}=1+\ba+\cA-\e{123}$.
The remaining scalar is  $a_{-}=2\sqrt{38}$, $(a_{0}-a_{123})=2$.
Then, Eq.~\eqref{B0ForCl03FullAplus} gives
$\m{A}_{{0}_{+}}=\log(0_{+})$; Eq.~\eqref{generic03allterms2}
gives $\m{A}_{{0}_{-}}=\tfrac{1}{2}\log 156$;
Eq.~\eqref{B12ForCl03FullPlus} gives $\m{A}_{{1,2}_{+}}=0$;
Eq.~\eqref{generic03allterms4} gives
$\m{A}_{{1,2}_{-}}=\frac{\arctan\left(\sqrt{38}\right)+2 \pi
c_{1_{-}}}{\sqrt{38}} (\ba+\cA)$. Finally, the logarithm of
$\m{A}$ is
\begin{equation}\label{ex3}
\begin{split}
  \log(\m{A})=&
  \frac{\arctan\left(\sqrt{38}\right)+2 \pi  c_{1_{-}}}{2\sqrt{38}} \left(\ba+\cA\right)\\
  &+ \tfrac12\big(\log(0_{+})\left(1+\e{123}\right)+\tfrac{1}{2}\log (156) \left(1-\e{123}\right)\big)\,.
\end{split}
\end{equation}
Note the factor $\log(0_{+})$ in front of
$\left(1+\e{123}\right)$. If logarithm in this form is inserted
into coordinate-free exponential~\cite{Acus2022} we will get
\begin{equation}
  \left(\tfrac{1}{2}\ee^{\log(0_{+})}+1\right)+\ba+\cA +\left(\tfrac{1}{2}
  \ee^{\log(0_{+})}-1\right)\e{123}\,,
\end{equation}
which coincides with the initial MV if we assume that\footnote{The
statement can be made strict by considering the limit $\lim_{x\to
0_{+}}\exp(\log(x))=0$, where $x\to 0_{+}$ indicates that the
limit is taken keeping $x$ positive, i.e."from above".}
$\log(0_{+})=-\infty$.

\subsection{GA Logarithm of blades and their combinations in \cl{0}{3}}\label{individualcl03}
In this subsection, the logarithms for individual blades and their
combinations that follow from generic logarithm
(Theorem~\ref{theoremCl03}),
and may be useful in practice are collected. The norms listed below  are positive scalars.\\
 \indent\textit{Vector norm}: $\magnitude{\ba}=\sqrt{\ba
\gradeinverse{\ba}}=\sqrt{a_{1}^2+a_{2}^2+a_{3}^2}$ \,.\\
 \indent\textit{Paravector norm}:
$\magnitude{a_0+\ba}=\magnitude{\m{A}_{0,1}} =\bigl(\m{A}_{0,1}
\gradeinverse{\m{A}}_{0,1}\bigr)^\frac{1}{2}=\sqrt{a_{0}^2+a_{1}^2+a_{2}^2+a_{3}^2}$\,.\\
 \indent\textit{Bivector norm}: $\magnitude{\cA}=\bigl(\cA
\reverse{\cA}\bigr)^\frac{1}{2}=\sqrt{a_{12}^2+a_{13}^2+a_{23}^2}$\,.\\
 \indent\textit{Rotor norm}: $\magnitude{a_0+\cA}=\magnitude{\m{A}_{0,2}}
=\bigl(\m{A}_{0,2}
\reverse{\m{A}}_{0,2}\bigr)^\frac{1}{2}=\sqrt{a_{0}^2+a_{12}^2+a_{13}^2+a_{23}^2}$\,.

\vspace{3mm}
 \noindent Logarithms of blades and their combinations.

 \textit{Logarithm of vector} $\ba=a_{1}\e{1}+a_{2} \e{2}+a_{3}
\e{3}$,  $c_i\in\bbZ$,
\begin{align}
  \log(\ba)=&
    \frac{1}{2}\log (\magnitude{\ba}^2)+ \pi\frac{\ba}{\magnitude{\ba}}\bigl(\frac{1}{2}+c_1(1+I)+c_2(1-I)\bigr)
   ,&\magnitude{\ba}^2\neq 0 .\label{VecForCl03}
\end{align}

\textit{Logarithm of paravector} $\m{A}_{0,1}=a_0+\ba$; $c_{i}\in
\mathbb{Z}$ and $\hat{\bu}^2=-1$, $\hat{\cU}^2=-1$.
\begin{align}\log\m{A}_{0,1}=&
  \begin{aligned}[t]&
  \frac{1}{2}\log (\magnitude{\m{A}_{0,1}}^2) +\frac{\ba}{\magnitude{\ba}}\bigl(\arctan(a_0,\magnitude{\ba})\\
    &\qquad\qquad  + \pi(c_1(1+I)+c_2(1-I))\bigr),
    \end{aligned}&\qquad \magnitude{\ba}\neq 0.
\end{align}

\textit{Logarithm of bivector} $\cA=a_{12} \e{12}+a_{13}
\e{13}+a_{23} \e{23}$, $c_i\in\bbZ$,
\begin{align}
  \log(\cA)=&
    \frac{1}{2}\log (\magnitude{\cA}^2)+
    \pi\frac{\cA}{\magnitude{\cA}}\bigl(\frac{1}{2}+c_1(1+I)+c_2(1-I)\bigr),&\magnitude{\cA}^2\neq 0.\label{BivForCl03}
\end{align}

\textit{Logarithm of parabivector and rotor} $\m{A}_{0,2}=a_0+\cA$,

\begin{align}\log\m{A}_{0,2}=&
  \begin{aligned}[t]&
  \frac{1}{2}\log (\magnitude{\m{A}_{0,2}}^2) +\frac{\cA}{\magnitude{\cA}}\bigl(\arctan(a_0,\magnitude{\cA})\\
    &\qquad \qquad + \pi(c_1(1+I)+c_2(1-I))\bigr),
    \end{aligned}&\qquad\magnitude{\cA}\neq 0.\label{ParaVecForCl03}
\end{align}

\textit{Logarithm of center} $\m{A}_{0,3}=a_0+a_{123} I$,
\begin{align}\log\m{A}_{0,3}=&\begin{cases}
  \begin{aligned}[t]&
    \bigl(\tfrac{1}{2}\log (a_0-a_{123})+\pi c_1 \hat{\cU}_1\bigr)\bigl(1-I\bigr)\\
    & +  \bigl(\tfrac{1}{2}\log (a_0+a_{123})+\pi c_2 \hat{\cU}_2\bigr)\bigl(1+I\bigr),\end{aligned}&
    \begin{aligned}[t](a_0-a_{123})&>0\textrm{ and }\\
    (a_0+a_{123})&>0\\\end{aligned}\\[17pt]
  \begin{aligned}[t]&
   \bigl(\tfrac{1}{2}\log (a_0-a_{123})+\pi c_1 \hat{\cU}_1\bigr)\bigl(1-I\bigr)\\
  & +  \bigl(\tfrac{1}{2}\log (-a_0-a_{123})+\pi (c_2+\tfrac{1}{2}) \hat{\bu}_2\bigr)\bigl(1+I\bigr),\end{aligned}
    &\begin{aligned}[t](a_0-a_{123})&>0\textrm{ and }\\ (a_0+a_{123})&<0\\\end{aligned}\\[17pt]
  \begin{aligned}[t]&
    \bigl(\tfrac{1}{2}\log (-a_0+a_{123})+\pi (c_1+\tfrac{1}{2}) \hat{\bu}_1\bigr)\bigl(1-I\bigr)\\
  & +  \bigl(\tfrac{1}{2}\log (a_0+a_{123})+\pi c_2 \hat{\cU}_2\bigr)\bigl(1+I\bigr)
    ,\end{aligned}&\begin{aligned}[t](a_0-a_{123})&<0\textrm{ and }\\ (a_0+a_{123})&>0\\\end{aligned}\\[17pt]
\begin{aligned}[t]&
    \bigl(\tfrac{1}{2}\log (-a_0+a_{123})+\pi (c_1+\tfrac{1}{2}) \hat{\bu}_1\bigr)\bigl(1-I\bigr)\\
  & +\bigl(\tfrac{1}{2}\log (-a_0-a_{123})+\pi (c_2+\tfrac{1}{2}) \hat{\bu}_2\bigr)\bigl(1+I\bigr)
    ,\end{aligned}&\begin{aligned}[t](a_0-a_{123})&<0\textrm{ and }\\ (a_0+a_{123})&<0\\\end{aligned}
\label{centerCL03}
\end{cases}
\end{align}
where $\hat{\bu}_i$ and $\hat{\cU}_j$ are arbitrary non-commuting
unit vector and bivector, respectively. If $(a_0-a_{123})=0$ or
$(a_0+a_{123})=0$ some of subparts give
  $\log(0_{+})$.

\section{MV  logarithms in \cl{3}{0} and \cl{1}{2}}
\label{sec:Cl30}
$\cl{3}{0}$ and
$\cl{1}{2}$ algebras are isomorphic. Their multiplication tables
  coincide, for example, after the following exchange of basis elements:
 \begin{equation*}\begin{split}
\cl{3}{0}\quad&\{1,\quad \e{1},\quad \e{2},\quad \e{3},\ \ \e{12},\ \ \e{13},\ \e{23},\ \e{123}\}\downarrow\\
\cl{1}{2}\quad&\{1,\ -\e{1},-\e{12},-\e{13},-\e{2},-\e{3},\
\e{23}, -\e{123}\}.
\end{split}\end{equation*}
To find  formulas for logarithm in coordinate form the same
inverse solution method was used as described in the
Appendix~\ref{coordFormCL03} for \cl{0}{3} algebra. The logarithm
in \cl{3}{0} and \cl{1}{2} exists for all MVs except of nonzero
MVs of the   form $\m{A}_{1,2}=\ba+\cA$ that  satisfy the
condition $\Det(\m{A}_{1,2})=(a_{+}^2+a_{-}^2)^2=0$, i.e., for MVs
that are the sums of vector and bivector and the  determinant are 
equal to zero. These restrictions are the same as those for
GA square root to  exist (see~\cite{AcusDargysPreprint2020} and
Example~3 herein in case $s=S=0$).

\subsection{Logarithm formula for generic MV}\label{allCasesCl30}
\begin{thm}[Logarithm of multivector in \cl{3}{0} and \cl{1}{2}]\label{theoremCl30}
The logarithm of generic MV $\m{A}=a_0+\ba+\cA+a_{123}I$ is
another MV
\begin{equation}
  \label{log30Generic}
  \log(\m{A})= \m{A}_{{0}} +\m{A}_{{1,2}_{\log}}+\m{A}_{{1,2}_{\arctan}}+\m{A}_{I},
\end{equation}
where
\begin{align}
\m{A}_{{0}}&=
\frac{1}{2} \bigl(\log k_{+}+\log k_{-}\bigr), & a_{+}^2 + a_{-}^2 \neq 0
 \label{log30GenericS}\allowdisplaybreaks\\
  \m{A}_{{1,2}_{{\log}}}&=
  \frac{1}{2} \frac{a_{+}-a_{-}I}{a_{-}^2+a_{+}^2} \bigl(\log k_{+}-\log k_{-}\bigr) \bigl(\ba+\cA\bigr), & a_{+}^2 + a_{-}^2\neq0
\label{log30GenericLog}
\end{align}
\begin{multline}
\m{A}_{{1,2}_{\arctan}}= I\frac{a_{+}-a_{-}I}{a_{-}^2+a_{+}^2}\bigl( \ba+\cA \bigr)
       \Bigl(
  \frac{1}{2}\arctan\bigl(-(a_{+}^2 - a_{0}^2) -  (a_{-}^2 - a_{123}^2), \\[-1pt]
\hfill         (a_{+} - a_0)(a_{-} +  a_{123}) - (a_{+} + a_0)(a_{-} - a_{123})\bigr)+2\pi c_1\Bigr),\\
  \hfill \textrm{when}\quad   (a_{+}^2 + a_{-}^2\neq0)\quad \textrm{and}\quad (k_{-}k_{+}\neq  0),\label{log30GenericArctan}
\end{multline}
\begin{multline}
  \m{A}_{I} =I\arctan \bigl((a_{+}+a_{0}) k_{-}-(a_{+}-a_{0}) k_{+},
     (a_{-}+a_{123}) k_{-}-(a_{-}-a_{123}) k_{+}\bigr)\\
    \qquad + 2\pi c_2 I,\qquad \hfil
  \textrm{when}\quad (a_{+}^2 + a_{-}^2\neq 0)\quad \textrm{and either}\\
  (a_{+}+a_{0}) k_{-}-(a_{+}-a_{0}) k_{+}\neq 0\quad \textrm{or}\quad (a_{-}+a_{123}) k_{-}-(a_{-}-a_{123}) k_{+}\neq 0 \label{log30GenericI}
\end{multline}
where scalar coefficients are
\begin{align}
  k_{-}^2&=(a_{+}-a_{0})^2+(a_{-}-a_{123})^2,\qquad
  k_{+}^2=(a_{+}+a_{0})^2+(a_{-}+a_{123})^2,\label{kdef}
\end{align}
and
\begin{align}
&a_{-}=\frac{-2 I \ba \w \cA}{\sqrt{2} \sqrt{\ba\d\ba+\cA\d\cA+\sqrt{(\ba\d\ba+\cA\d\cA)^2 -4  (\ba\w\cA)^2}}},\notag\\
&a_{+}=\frac{\sqrt{\ba\d\ba+\cA\d\cA+\sqrt{(\ba\d\ba+\cA\d\cA)^2 -4  (\ba\w\cA)^2}}}{\sqrt{2}}\,\label{bPM} \\
&\quad \textrm{for } \ba \w \cA \neq 0\notag, \quad \textrm{and}\allowdisplaybreaks\\
&\begin{cases}
  a_{+}=\sqrt{\ba\d\ba+\cA\d\cA},\quad a_{-}=0,&\textrm{if}\quad \ba\d\ba+\cA\d\cA \ge 0\\
a_{+}=0,\quad a_{-}=\sqrt{-(\ba\d\ba+\cA\d\cA)},&\textrm{if}\quad
\ba\d\ba+\cA\d\cA < 0,\notag
\end{cases}\\
 &\qquad \textrm{when}\quad \ba \w \cA = 0\notag.
\end{align}
  The constants $c_{1},c_{2}$ are arbitrary integers.
\end{thm}
\textit{Proof}. It is enough to check that after substitution of
$\log{\m{A}}$ expressions into exponential formula presented
in~\cite{Dargys2022a} one gets the initial MV~$\m{A}$. The
factor $\frac{a_{+}-a_{-}I}{a_{-}^2+a_{+}^2}$ in the above
formulas alternatively may be written as $(a_{+}+a_{-}I)^{-1}$.

\subsection{Special cases}

When the conditions listed in
Eqs~\eqref{log30GenericS}-\eqref{log30GenericI} are not satisfied,
we have  special cases. In particular, the condition $k_{\pm}=0$
means that the MV determinant  is zero, $\Det(\m{A})= k_{-}^2
k_{+}^2=0$. Similarly, the condition $a_{+}^2 + a_{-}^2=0$ implies
that determinant of vector+bivector part  vanishes,
$\Det(\m{A}_{1,2})=(a_{+}^2 + a_{-}^2)^2=0$. The specific
relations $(a_{+}+a_{0}) k_{-}-(a_{+}-a_{0}) k_{+}\neq 0$ and
$(a_{-}+a_{123}) k_{-}-(a_{-}-a_{123}) k_{+}\neq 0$ in
Eq.~\eqref{log30GenericI} as well as the relation $k_{-}k_{+}\neq
0$ in Eq.~\eqref{log30GenericArctan} ensure that both arguments of
$\arctan(x,y)$ do not nullify simultaneously.

When the generic formula is not applicable the expressions for
$\m{A}_{{0}}, \m{A}_{{1,2}_{{\log}}}, \m{A}_{{1,2}_{\arctan}}$ and
$\m{A}_{I}$ in Theorem~\ref{theoremCl30} must be supplemented by
following formulas

\begin{align}
\begin{split}\label{log30AllS}
  \m{A}_{{0}}{}&=\begin{cases}
    \frac{1}{2} \log\bigl(a_0^2+a_{123}^2\bigr), &  \begin{aligned}[t](a_{+}^2 + a_{-}^2 =& 0)\land (a_{0}^2 + a_{123}^2\neq 0),\end{aligned}\\
    \varnothing,& \begin{aligned}[t] (a_{+}^2 + a_{-}^2 = 0)\land (a_{0}^2 + a_{123}^2= 0)\end{aligned}
\end{cases}
\end{split}\\
\begin{split}\label{log30AllLog}
  \m{A}_{{1,2}_{{\log}}}{}&=\begin{cases}
  0, & \begin{aligned}[t](a_{+}^2 + a_{-}^2 = 0)\land
(a_{0}^2 + a_{123}^2\neq 0),\end{aligned}\\
   \varnothing,& \begin{aligned}[t] (a_{+}^2 + a_{-}^2 = 0)\land (a_{0}^2 + a_{123}^2= 0)\end{aligned}
\end{cases}
\end{split}
\\
\begin{split}\label{log30AllArctan}
  \m{A}_{{1,2}_{\arctan}}{}&=
  \begin{cases}
  \pi(\frac{1}{2}+2c_1)I\frac{a_{+}-a_{-}I}{a_{-}^2+a_{+}^2}\bigl( \ba+\cA \bigr), & \begin{aligned}[t]
  (a_{+}^2 + a_{-}^2\neq0)\land(k_{-}k_{+}=  0),\end{aligned}\\
    \frac{a_{0}-a_{123}I}{a_{0}^2+a_{123}^2}\bigl( \ba+\cA \bigr)+\hat{\cF}, & \begin{aligned}[t]&(a_{+}^2 + a_{-}^2 = 0)\land (a_{0}^2 + a_{123}^2\neq 0),\end{aligned}\\
    \varnothing,& \begin{aligned}[t](a_{+}^2 + a_{-}^2 = 0)&\land (a_{0}^2 + a_{123}^2= 0)\end{aligned}
\end{cases}
\end{split}\\
\begin{split}\label{log30AllI}
  \m{A}_{I}{}&=\begin{cases}
  I\bigl( \arctan(-a_{-},a_{+})+2\pi c_2\bigr), &\begin{aligned}[t]
    &(a_{+}^2 + a_{-}^2\neq0)\\
    &\land
    ((a_{+}+a_{0}) k_{-}-(a_{+}-a_{0}) k_{+} = 0)\\
    &\land((a_{-}+a_{123}) k_{-}-(a_{-}-a_{123}) k_{+} = 0),
 \end{aligned}\\
    I \bigl(\arctan(a_{0},a_{123})+2\pi c_2\bigr), & \begin{aligned}[t]& (a_{+}^2 + a_{-}^2 = 0)\land(a_{0}^2 + a_{123}^2\neq 0),\end{aligned}
      \\
\varnothing,& \begin{aligned}[t] (a_{+}^2 + a_{-}^2 = 0)&\land (a_{0}^2 + a_{123}^2= 0)\end{aligned}
\end{cases}
\end{split}
\end{align}
Here the symbols $\land$ and $\lor$ in the conditions represent
logical conjunction and disjunction, respectively. $\hat{\cF}=
\begin{cases} 2\pi c_{1}\hat{\cU},& \mathrm{if}\quad \ba+\cA = 0 \\
    0, &\mathrm{if}\quad \ba+\cA\neq 0
  \end{cases}$, where the free unit bivector 
must satisfy $\hat{\cU}^2=-1$. After exponentiation it gives
$\exp(\hat{\cU})=1$ and represents  continuous degree of freedom
(a direction) in \eqref{log30AllArctan} and \eqref{log30AllI}, and
can be parameterized as
\begin{equation}
  \hat{\cU}=\begin{cases}\frac{d_{12}\e{12}+d_{13}\e{13}+d_{23}\e{23}}{\sqrt{d_{12}^2+d_{13}^2+d_{23}^2}},& \mathrm{for}\quad \cl{3}{0}\\
\frac{d_{12}\e{12}+d_{13}\e{13}+d_{23}\e{23}}{\sqrt{-d_{12}^2-d_{13}^2+d_{23}^2}},&
\mathrm{for}\quad \cl{1}{2},\ \mathrm{when}\
-d_{12}^2-d_{13}^2+d_{23}^2 > 0
\end{cases}
\end{equation}
The cases $k_{\pm}=0$ that represent MV with a vanishing
determinant, $\Det(\m{A})= k_{-}^2  k_{+}^2=0$,  yield MVs with
infinite coefficients (see Example~9 for details).

\subsection{Multiveluedness and free multivector}\label{freeParametersCl30}

In Eqs~\eqref{log30GenericArctan} and \eqref{log30GenericI}  we
may add any multiple of $2\pi$ to both  arc tangent functions,
i.e. $\arctan (y_1,y_2)\to \arctan (y_1,y_2)+2\pi c_i$. After
collecting terms in front of free coefficients $c_1,c_2\in
\mathbb{Z}$, we obtain a free MV $\m{F}$ that satisfies
$\exp(\m{F})=1$,
\begin{align}\label{periodicityCL30CL12}
  \m{F}=& \frac{2\pi c_1}{\bigl(a_{-}^2+a_{+}^2\bigr)} \big(a_{-}
  (\ba+\cA)+a_{+} (\ba+\cA) I\big) +2\pi c_2 I\, ,
\end{align}
where $a_{\pm}$ are given by Eq.~\eqref{bPM}.

\vspace{3mm} \textbf{Example 7.} {\it Logarithm of generic MV in
\cl{3}{0}.} Let us take simple but representative MV:
$\m{A}=-2+\e{1}+ \e{23}-3 \e{123}$. From Eqs~\eqref{kdef} and
\eqref{bPM} we have $k_{+}^2=5$, $k_{-}^2=25$ and $a_{+}=a_{-}=1$.
Then~\eqref{log30GenericS} and \eqref{log30GenericLog} yield
$\m{A}_{{0}}=\frac{3\log 5}{4}$ and
$\m{A}_{{1,2}_{\log}}=-\frac{\log 5}{8} \bigl(\e{1}+
\e{23}\bigr)(1-I)$. Next, the Eqs~\eqref{log30GenericArctan} and
\eqref{log30GenericI} give
$\m{A}_{{1,2}_{\arctan}}=-\frac{1}{4}\bigl(-\arctan\frac{2}{11}+4\pi
c_2\bigr) \bigl(\e{1}+ \e{23}\bigr)(1+I)$ and
$\m{A}_{I}=\Bigl(-\pi+\arctan\frac{-10-4 \sqrt{5}}{-5-3 \sqrt{5}}
+2\pi c_1\Bigr)\e{123}$. Finally, after summation of all terms
in~\eqref{log30Generic} we obtain $\log(\m{A})=\frac{\log
5}{4}\bigl(3-\e{1}\bigr)+\frac{1}{2}\arctan\frac{2}{11}\e{23}
+\bigl(-\pi+\arctan(\frac{1}{2}(1+\sqrt{5})\bigr)\e{123}+\m{F}$,
where the free MV $\m{F}=2\pi\bigl(c_1\e{123}-c_2\e{23}\bigr)$.
The coefficients $c_1,c_2\in \mathbb{Z}$ come from
$\m{A}_{{1,2}_{\arctan}}$ and $\m{A}_{I}$ terms, respectively.
Substitution of this result into exponential $\exp(\log(\m{A}))$
returns the initial MV.

\textbf{Example 8.} {\it Logarithm of center of \cl{3}{0}.}
$\m{A}=1- 2 \e{123}$. Since $\e{123}^2=-1$ the MV is a
counterpart of complex number logarithm.  Eqs~\eqref{kdef}
and~\eqref{bPM} give $a_{+}=a_{-}=0$ and $k_{+}^2=k_{-}^2=5$.
Then, Eq.~\eqref{log30AllS} gives $\m{A}_{{0}}=\frac{\log 5}{2}$;
Eq.~\eqref{log30AllLog} gives $\m{A}_{{1,2}_{\log}}=0$;
Eq.~\eqref{log30AllArctan} gives $\m{A}_{{1,2}_{\arctan}}=2\pi
c_1\hat{\cU}$; Eq.~\eqref{log30AllI} gives
$\m{A}_{I}=\bigl(-\arctan 2 +2\pi c_2\bigr)\e{123}$. Note that
$\hat{\cU}$ is the same free MV for both
$\m{A}_{{1,2}_{\arctan}}$. After summation of  terms
in~\eqref{log30Generic} the final answer is
$\log(\m{A})=\frac{\log 5}{2}+\bigl(-\arctan 2 +2\pi
c_2\bigr)\e{123}+2\pi c_1\hat{\cU}$. On the other hand the
complex number $1-2\,\ii$ gives $\log(1-2\,\ii)=\tfrac12\log
5-\arctan 2$, which coincides with  \cl{3}{0} algebra result if
$c_1=c_2=0$.

\textbf{Example 9.} {\it Logarithm of singular MV when
$\Det(\m{A})=0$.} This is the most intriguing and complicated case
in \cl{3}{0}. Since  $\Det(\m{A})= k_{-}^2  k_{+}^2$ we may have
either $k_{-}^2=0$ or $k_{+}^2=0$. The case when
$k_{-}^2=k_{+}^2=0$ is trivial since it requires all MV components
to vanish. Let's analyze the case when $k_{+}^2\neq0$ and
$k_{-}^2=0$. It is represented, for example, by
$\m{A}=6+(-8\e{1}-2\e{3})+(-\e{12}+ 10\e{13}+ 10\e{23})-13
\e{123}=6+\ba+\cA-13 \e{123}$. From Eq.~\eqref{bPM} we find
$a_{+}=6$, $a_{-}=-13$ and from~\eqref{kdef} $k_{+}^2=820$,
$k_{-}^2=0$. Then, Eq.~\eqref{log30GenericS} gives
$\m{A}_{{0}}=\frac{1}{2} \bigl(\log(2\sqrt{205}) + \log
(0_{+})\bigr)$; Eq.~\eqref{log30GenericLog} gives
$\m{A}_{{1,2}_{{\log}}}=\frac{1}{410} \bigl(\log(2\sqrt{205}) -
\log (0_{+})\bigr) \bigl(6+13\e{123}\bigr)\bigl(\ba
+\cA\bigr)+6(\ba+\cA) \bigr)$; Eq.~\eqref{log30AllArctan} gives
$\m{A}_{{1,2}_{\arctan}}= \frac{\pi}{205}(\frac{1}{2}+2
c_1)\bigl(-6+13\e{123}\bigr)\bigl(\ba+\cA\bigr)+6(\ba +\cA)
\bigr)$; finally, Eq.~\eqref{log30AllI} gives
 $\m{A}_{I}=\bigl(\arctan(\frac{6}{13})+2\pi c_2\bigr)\e{123}$.
Summing up all terms we obtain the answer
\begin{align*}
\log \m{A}=&\frac{1}{2} \bigl(\log(2\sqrt{205}) + \log (0_{+})\bigr)
  + \Bigl(\frac{1}{410} \bigl(\log(2\sqrt{205}) - \log (0_{+})\bigr) \bigl(6+13\e{123}\bigr)\\
  &\quad + \frac{\pi}{205}(\frac{1}{2}+2 c_1)\bigl(-6+13\e{123}\bigr) \Bigr)(\ba+\cA)
  +\bigl(\arctan(\frac{6}{13})+2\pi c_2\bigr)\e{123}.
 \end{align*}
The result can be checked after replacement of  $\log(0_{+})$ by
$\log (x)$ and substitution into exponential formula (4.1) of
paper~\cite{Dargys2022a}. After simplification one can take the
limit $\lim_{x\to 0_{+}}\exp\bigl(\log \m{A}\bigr)$, which returns
the initial~MV. This example demonstrates that the logarithm of MV
with specific finite coefficients may yield MV with some of
coefficients in the answer being infinite and which have to be
understood as the limit $\lim_{x\to 0_{+}}\log(x)$. The answer,
nevertheless, is meaningful since the substitution of the answer
into exponential formula and computation  of the limit reproduces
the initial MV.

\subsection{Logarithms of individual blades and their combinations}\label{individualcl30}

Below we use different norms for individual blades of
\cl{3}{0}, since a positive scalar for vectors and bivectors is
calculated differently. In particular, for a vector we will use
$\magnitude{\ba}=\sqrt{\ba \ba}=\sqrt{a_{1}^2+a_{2}^2+a_{3}^2}$\,,
whereas for a bivector $\magnitude{\cA}=\bigl(\cA
\reverse{\cA}\bigr)^{1/2}=\sqrt{a_{12}^2+a_{13}^2+a_{23}^2}$\,.
For a rotor $\magnitude{a_0+\cA}=\magnitude{\m{A}_{0,2}}
=\bigl(\m{A}_{0,2}
\reverse{\m{A}}_{0,2}\bigr)^{1/2}=\sqrt{a_{0}^2+a_{12}^2+a_{13}^2+a_{23}^2}$\,,
and for an element of center
$\magnitude{a_0+a_{123}I}=\magnitude{\m{A}_{0,3}}=\bigl(\m{A}_{0,3}
\gradeinverse{\m{A}}_{0,3}\bigr)^{1/2}=\sqrt{a_{0}^2+a_{123}^2}$\,.

Logarithm of {\it vector}: $\ba=a_{1}\e{1}+a_{2} \e{2}+a_{3}
\e{3}$,
\begin{align}
  \log(\ba)=&
    \frac{1}{2}\log (\magnitude{\ba}^2)-
    \pi\bigl(\tfrac{1}{2}+2 c_2\bigr)\frac{\ba}{\magnitude{\ba}}I+\pi\bigl(\tfrac{1}{2}+2 c_1\bigr) I,&\magnitude{\ba}^2\neq 0\,.\label{VecForCl30}
\end{align}

Logarithm of {\it bivector}: $\cA=a_{12} \e{12}+a_{13}
\e{13}+a_{23} \e{23}$,
\begin{align}
  \log(\cA)=&
    \frac{1}{2}\log (\magnitude{\cA}^2)-
    \pi\bigl(\tfrac{1}{2}+2 c_2\bigr)\frac{\cA}{\magnitude{\cA}}+\pi(1+2 c_1) I,&\magnitude{\cA}^2\neq 0\,.\label{BivForCl30}
\end{align}

Logarithm of \textit{rotor}: $\m{A}_{0,2}=a_0+\cA$,
\begin{align}\log\m{A}_{0,2}=&\begin{cases}
  \begin{aligned}&
  \frac{1}{2}\log (\magnitude{\m{A}_{0,2}}^2) +\bigl(\arctan\bigl(a_0,0\bigr)+2\pi c_1\bigr)I\\
    &\quad   +\frac{\cA}{\magnitude{\cA}}\Bigl(2\pi c_2-\frac{1}{2}\arctan\bigl(a_0^2-\magnitude{\cA}^2,-2a_0\magnitude{\cA}\bigr)\Bigr)
    \end{aligned},
    &\begin{aligned}\magnitude{\m{A}_{0,2}}\neq 0\end{aligned}\\
      \log a_0+2\pi c_1 I,&\magnitude{\cA}= 0\textrm{ and } a_{0}\ge 0\\
      \log (-a_0)+2\pi (c_1+1) I,&\magnitude{\cA}= 0\textrm{ and } a_{0}< 0\\
  \text{see bivector formula~\eqref{BivForCl30}},&\magnitude{\cA}\neq 0\textrm{ and }a_{0}=0,
  \label{RotorCl30}
\end{cases}
\end{align}

Logarithm of  \textit{center}: $\m{A}_{0,3}=a_{0}+a_{123}
\e{123}=a_{0}+a_{123} I$,
$\magnitude{\m{A}_{0,3}}^2=\m{A}_{0,3}\widetilde{\m{A}}_{0,3}$.
\begin{align}
  \log(\m{A}_{0,3})=&\begin{cases}
    \frac{1}{2}\log (\magnitude{\m{A}_{0,3}}^2)+ 2\pi c_2\hat{\cU}+
     \bigr(\arctan\bigl(a_0,a_{123}\bigr)+4\pi c_1\bigl)I,&\magnitude{\m{A}_{0,3}}^2\neq 0,\\
    \log(0_{+})+ 2\pi c_2\hat{\cU},& \magnitude{\m{A}_{0,3}}^2 = 0\,.\label{CenterForCl30}\\
\end{cases}
\end{align}
The paravector $\m{A}_{0,1}=a_0+\ba$ norm
$\magnitude{a_0+\ba}^2\equiv\magnitude{\m{A}_{0,1}}^2
=\m{A}_{0,1}\gradeinverse{\m{A}}_{0,1}=a_{0}^2-a_{1}^2-a_{2}^2-a_{3}^2$,
contains coefficients with opposite signs. The logarithm formula,
therefore, splits into many subcases and is impractical.

\section{MV  logarithms in \cl{2}{1}}
\label{sec:Cl21}

Of all three algebras, the logarithm of \cl{2}{1} appeared the
most hard to recover.  The logarithms in $\cl{3}{0}$ and
$\cl{1}{2}$ algebras exist for almost all MVs except very small
specific class of vectors and bivectors, $\ba+\cA\neq 0$, with the
vanishing determinant $\Det{(\ba+\cA)}= 0$. In $\cl{2}{1}$ algebra
the logarithm does not exist for a large class of MVs. In
contrast, in $\cl{0}{3}$ algebra the logarithm exists for all MVs.
\begin{thm}\label{theorem3}[Logarithm of multivector in \cl{2}{1}]
  \label{log21thmTheorem} The logarithm of multivector $\m{A}=a_0+(a_1\e{1}+a_2\e{2}+a_3\e{3})+(a_{12}\e{12}+a_{13}\e{13}+a_{23}\e{23})+a_{123}I=a_0+\ba+\cA+a_{123}I$ is the MV
  \enlargethispage{5pt}
 \begin{align}
    \log(\m{A})= &\begin{cases}\frac{1}{2}\bigl(\m{A}_{{0}_{+}}+\m{A}_{{0}_{-}}+\m{A}_{{1,2}_{+}}+\m{A}_{{1,2}_{-}}+(\m{A}_{{0}_{+}} -\m{A}_{{0}_{-}})I\bigr),& f_{\pm}\ge 0 \\
    \varnothing,& f_{\pm}< 0\end{cases}\label{LogCl21Full}\allowdisplaybreaks\\
\textrm{where}\notag
\\
&\begin{aligned}
  &f_{\pm}=(a_{0}\pm a_{123})^2+a_{\pm}^2,&f_{\pm}\lesseqqgtr 0,\\
  &a_{-}^{(2)}=-(\ba\d\ba+\cA\d\cA) + 2I  \ba\w\cA,&a_{-}^{(2)}\lesseqqgtr 0,\label{apliusminuscl21}\\
  &a_{+}^{(2)}=-(\ba\d\ba+\cA\d\cA) - 2I  \ba\w\cA,&a_{+}^{(2)}\lesseqqgtr 0,\allowdisplaybreaks\\
\end{aligned}\\
\textrm{and}\notag\\
   \m{A}_{{0}_{\pm}}=&\begin{cases}
    \tfrac{1}{2}\log(f_{\pm}),&\begin{aligned}[t] (a_{\pm}^{(2)} > 0)\end{aligned}\\
      \begin{aligned}[t]&\tfrac{1}{2}\log\Bigl(a_0\pm a_{123}+\sqrt{-a_{\smash{\pm}}^{(2)}}\Bigr)\\ & + \tfrac{1}{2}\log\Bigl(a_0\pm a_{123}-\sqrt{-a_{\smash{\pm}}^{(2)}}\Bigr),\end{aligned}
    &\begin{aligned}[t](a_{\pm}^{(2)} < 0)\land (a_{0}\pm a_{123} > 0)\end{aligned}\\
      \log(a_{0}\pm a_{123})+  2\pi c_{2\pm}\hat{\cF},& (a_{\pm}^{(2)}=0)\land (a_{0}\pm a_{123} > 0)\\
      \log\bigl(-(a_{0}\pm a_{123})\bigr)+  (\pi+2\pi c_{2\pm})\hat{\cU},& \begin{aligned}[t](a_{\pm}^{(2)}=0)&\land (a_{0}\pm a_{123} \le 0)\\& \land(\mathfrak{D}= \mathrm{True})\end{aligned}\\
    \varnothing,&\begin{aligned}[t]&\bigl((a_{\pm}^{(2)} <  0) \land (a_{0} \pm a_{123} < 0)\bigr)\\
      &\lor\bigl( (a_{\pm}^{(2)} =  0) \land (a_{0} \pm a_{123} \le 0) \\
    &\hphantom{\lor \bigl( (a_{\pm}^{(2)} =  0)\land}\land (\mathfrak{D}= \mathrm{False})\bigr)\end{aligned}
    \label{B0ForCl21FullFnonzero}\\
\end{cases}\allowdisplaybreaks\\
  \m{A}_{{1,2}_{\pm}}=&\begin{cases}
    \frac{1}{\sqrt{a_{\smash{\pm}}^{(2)}}}\bigl(\arctan (a_{0}\pm a_{123},a_{\pm})\\
    \qquad  +2\pi c_{1\pm} \bigr) (1\pm I)\bigl(\ba +\cA\bigr),&\begin{aligned}[t]&(a_{\pm}^{(2)}> 0)\end{aligned}\\
\begin{aligned}[t]&
  \frac{1}{\sqrt{-a_{\smash{\pm}}^{(2)}}}\arctanh \bigl(\frac{\sqrt{-a_{\smash{\pm}}^{(2)}}}{a_{0}\pm a_{123}}\bigr)\\
&\qquad\qquad\times (1\pm I)\bigl(\ba +\cA\bigr)\end{aligned},&\begin{aligned}[t](a_{\pm}^{(2)}< 0)&\land (a_{0}\pm a_{123} > 0)\\
& \land\bigl(-a_{\pm}^{(2)}\neq (a_{0}\pm a_{123})\bigr)\end{aligned}\\
      \begin{aligned}[t]&\tfrac{1}{2}\Bigl(\log\bigl(a_0\pm a_{123}+\sqrt{-a_{\smash{\pm}}^{(2)}}\bigr)\\ &\  - \log\bigl(a_0\pm a_{123}-\sqrt{-a_{\smash{\pm}}^{(2)}}\bigr)\Bigr)\\
   &
    \quad\times\frac{1}{\sqrt{-a_{\smash{\pm}}^{(2)}}}(1\pm I)\bigl(\ba +\cA\bigr)
    \end{aligned},&\begin{aligned}[t](a_{\pm}^{(2)}< 0)&\land (a_{0}\pm a_{123} > 0)\\&\land \bigl(-a_{\pm}^{(2)}= (a_{0}\pm a_{123})\bigr)\end{aligned}\\
     \frac{1}{a_{0}\pm a_{123}}  (1\pm I)\bigl(\ba +\cA\bigr),& (a_{\pm}^{(2)}=0)\land (a_{0}\pm a_{123} > 0)\\[2pt]
    0,& \begin{aligned}[t](a_{\pm}^{(2)}=0)&\land (a_{0}\pm a_{123} \le 0)\\&\land (\mathfrak{D}= \mathrm{True})\end{aligned}\\[2pt]
 \varnothing,&\begin{aligned}[t]&\bigl((a_{\pm}^{(2)} <  0) \land (a_{0} \pm a_{123} < 0)\bigr)\\
   &\lor \bigl( (a_{\pm}^{(2)} =  0) \land (a_{0} \pm a_{123} \le 0) \\
 &\hphantom{\lor \bigl((a_{\pm}^{(2)} <  0)}\land (\mathfrak{D}= \mathrm{False})\bigr)\end{aligned}
    \label{B12ForCl21FullFnonzero} \\
\end{cases}
\end{align}
    where the upper symbol in $a_{\pm}^{(2)}$ indicates that $a_{\pm}^{(2)}$ consists of the
  squared coefficients  $a_i^2$, $a_{ij}^2$ and $a_ia_{ij}$.
  $\hat{\cF}=\begin{cases} \hat{\cU},& \mathrm{if}\quad\mathfrak{D} = \mathrm{True} \\
    0, &\mathrm{if}\quad\mathfrak{D} = \mathrm{False}
  \end{cases}$. The logical condition $\mathfrak{D}$ is a conjunction of outcomes of three
  comparisons
    $\mathfrak{D}=(a_{1}=\pm a_{23})\land (a_{2}=\mp a_{13})\land (a_{3}=\mp a_{12})\equiv\bigl((a_{1}= a_{23})\land (a_{2}=- a_{13})\land (a_{3}=- a_{12})\bigr)\lor\bigl((a_{1}=- a_{23})\land (a_{2}= a_{13})\land (a_{3}= a_{12})\bigr)$
    that should be applied to $\m{A}_{{0}_{\pm}}$ and $\m{A}_{{1,2}_{\pm}}$ terms without paying attention to $\pm$ signs in their subscripts. Unit bivector in $\m{A}_{{0}_{\pm}}$ may be parameterized as
  $\hat{\cU}=\frac{d_{12}\e{12}+d_{13}\e{13}+d_{23}\e{23}}{\sqrt{d_{12}^2-d_{13}^2-d_{23}^2}}$.
    The symbol $\varnothing$ means that the solution set is empty.
    In all formulas the indices and conditions (except $\mathfrak{D}$ as  stated explicitly)  must be included with either all upper or with all lower signs.
\end{thm} 
The case with $f_{\pm}\neq 0$ and $a_{\pm}^{(2)} > 0$ represents a generic instance.
 When either $f_{\pm} = 0$ or $a_{\pm}^{(2)}\le 0$
we have the special case. Note  that in
Eq.~\eqref{apliusminuscl21} the condition $f_{\pm}=0$  implies
$a_{\pm}^{(2)} \le 0$. Also, observe that the condition
$f_{\pm}\ge 0$  ensures  automatically that a less restrictive
requirement $\Det(\m{A})=f_{-}f_{+}\ge0$ is  fulfilled
automatically.

The equations
\eqref{B0ForCl21FullFnonzero}-\eqref{B12ForCl21FullFnonzero} are
similar to
Eqs~\eqref{B0ForCl03FullAplus}-\eqref{B12ForCl03FullMinus}  in
\cl{0}{3} (see Sec.~\ref{spacialCL04}). Also,
in~\eqref{apliusminuscl21} the expressions for scalar
coefficients $a_{\pm}=\begin{cases}\sqrt{a_{\pm}^{(2)}},&a_{\pm}^{(2)}\ge 0\\
  \sqrt{-a_{\pm}^{(2)}},&a_{\pm}^{(2)}<0
\end{cases}$ are similar  to
Eqs~\eqref{aminusCl03} and \eqref{apliusCl03}. The differences
mainly arise  at the parameter boundaries that define the
existence of MV logarithm for $\cl{2}{1}$.

From our earlier calculations~\cite{AcusDargysPreprint2020} we
know the conditions that ensure an existence of MV square roots in
$\cl{2}{1}$ algebra. Thus, we can rewrite and use here these
conditions  that limit the extent of the logarithm in
Theorem~\ref{theorem3}. It appears that the quantities $b_S$ and
$b_I$ in~\cite{AcusDargysPreprint2020} may be expressed in terms
of multipliers $f_{+}$ and $f_{-}$ in the determinant $D=\Det
\m{A}=f_{-}f_{+}$, where $f_{\pm}=(a_{0}\pm
a_{123})^2+a_{\pm}^{(2)}$, in a form
$b_I=\frac{1}{2}\bigl(f_{+}-f_{-}\bigr)$ and
$b_S=\frac{1}{2}\bigl(f_{+}+f_{-}\bigr)$. Now, note that $f_{\pm}$
enter as arguments in log-functions of
Theorem~\ref{theorem3}, Eq.~\eqref{B12ForCl21FullFnonzero}.
Therefore, the square root existence condition $b_S-\sqrt{D}\ge 0$
in~\cite{AcusDargysPreprint2020}, in terms of the logarithm
problem can be rewritten as a difference of the determinant
factors, namely, $b_S-\sqrt{D}\Leftrightarrow
\frac{1}{2}\bigl(\sqrt{f_{-}}-\sqrt{f_{+}}\bigr)^2$. Now it
becomes clear that this condition is always satisfied and
therefore can be ignored, once we assume that the both factors
satisfy $f_{-}\ge 0$ and $f_{+}\ge 0$. From all this we conclude
that the requirement $f_{\pm}> 0$ constitutes one of the existence
conditions of logarithm in Theorem~\ref{theorem3}. Also,
$b_S-\sqrt{D}=0$ is equivalent to $f_{+}=f_{-}$. This restricts
the maximal possible value of $a_{\pm}^{(2)}$. In particular,
$|a_{\pm}^{(2)}|\le (a_{0}\pm a_{123})^2$. Remember, that notation
$a_{\pm}^{(2)}$ (instead of $a_{\pm}$) was introduced to keep an
analogy with $\cl{0}{3}$ case. It may be negative $a_{\pm}^{(2)} <
0$ (see definition~\eqref{apliusminuscl21}) and therefore the
notation, in general, can't be interpreted as a square of scalar
unless $a_{\pm}^{(2)} \ge 0$. When $a_{\pm}^{(2)}= 0$, an
additional condition $a_{0}\pm a_{123} \ge 0$ is required for
logarithm to exist.

Since \cl{2}{1} algebra is rarely used we will not provide
explicit formulas for pure blades (they can be found in the
notebook {\tt ElementaryFunctions.nb} in \cite{AcusDargys2017}).
Also, because generic formulas are similar to those in $\cl{0}{3}$
the examples below are restricted to special cases only.

\vspace{3mm} \textbf{Example 10.} {\it Logarithm in \cl{2}{1} when
$a_{\pm}^{(2)}=0$ and $a_0\pm a_{123}>0$.}  Let the MV be
$\m{A}=7+(2\e{1}+\e{2}+3\e{3})+(2\e{12}+ 2\e{13}-2\e{23})+5
I=7+\ba+\cA+5 I$. From~\eqref{apliusminuscl21} we find
$a_{+}^{(2)}=0,f_{+}=144$ and $a_{-}^{(2)}=0,f_{-}=144$. Since
$a_0\pm a_{123}=7\pm5>0$ from~\eqref{B0ForCl21FullFnonzero} we
have $\m{A}_{{0}_{-}}=\log 2$,  $\m{A}_{{0}_{+}}=\log 12$ and
from~\eqref{B12ForCl21FullFnonzero} $\m{A}_{{1,2}_{-}}=\frac12
(1-I)\bigl(\ba+\cA\bigr)$, $\m{A}_{{1,2}_{+}}=\frac{1}{12}
(1+I)\bigl(\ba+\cA\bigr)$. Finally, $\log \m{A}=
\frac{1}{24}\bigl(12\log(24)+24\e{1}+17\e{2}+31\e{3}+29\e{12}+
19\e{13}-24\e{23}+12\log(6) I\bigr)$. Note, because MV
coefficients $a_3\ne\pm a_{12}$  the condition $\mathfrak{D}$ is
$\mathrm{False}$, therefore the free MV in~
\eqref{B0ForCl21FullFnonzero} is absent, $\hat{\cF}=0$.

\textbf{Example 11.} {\it Logarithm  when $a_{-}^{(2)}=0$, $a_0-
a_{123}=0$ and $a_{+}^{(2)}>0$, $a_0+ a_{123}<0$.} In  \cl{2}{1}
these properties are satisfied by MV
$\m{A}=-2+(7\e{1}+4\e{2}+10\e{3})+(-10\e{12}-4\e{13}+7\e{23})-2
I=-2+\ba+\cA-2I$. From~\eqref{apliusminuscl21} we find
$a_{+}^{(2)}=140,f_{+}=156$ and $a_{-}^{(2)}=0,f_{-}=0$. Then,
because $a_0- a_{123}=-2-(-2)=0$ and $a_0+ a_{123}=-2-2<0$,
from~\eqref{B0ForCl21FullFnonzero} we have $\m{A}_{{0}_{-}}=\log
(0_{+})+(\pi+2\pi c_{2_{-}})\hat{\cU},
\m{A}_{{0}_{+}}=\frac12\log(156)$.
From~\eqref{B12ForCl21FullFnonzero} $\m{A}_{{1,2}_{-}}=0,
\m{A}_{{1,2}_{+}}=\frac{2}{\sqrt{35}}\bigl(\pi + 2\pi c_{1_{+}} -
\arctan(\sqrt{35}/2)\bigr) (1+I)\bigl(\ba+\cA\bigr)$. Then  using
\eqref{LogCl21Full} we obtain the final answer
$\log\m{A}=\alpha_0+\alpha_1
\e{1}+\alpha_2\e{2}+\alpha_3\e{3}+\alpha_{12}
\e{12}+\alpha_{13}\e{13}+\alpha_{23}\e{23}+\alpha_{123}I$, where
$\beta=\arctan(\sqrt{35}/2)$ and,
\[\begin{split}
&\alpha_0=\tfrac12(\log0_{+}+\log\sqrt{156}),\quad \alpha_{123}=-\tfrac14(2\log0_{+}-\log156),\\
&\alpha_1=-\tfrac{1}{20}\big(5\sqrt{3}\,\pi-3\sqrt{35}\,\pi+2\sqrt{35}\,\beta\big),\quad
\alpha_2=\big(\tfrac{\pi}{2\sqrt{3}}+\tfrac{2}{\sqrt{35}}\,\pi-\tfrac{2}{\sqrt{35}}\beta\big),\\
&\alpha_3=\big(\sqrt{\tfrac{5}{7}}\,\pi+\tfrac{5\pi}{4\sqrt{3}}-\sqrt{\tfrac{5}{7}}\,\beta\big),
\quad\alpha_{12}=\big(-\sqrt{\tfrac{5}{7}}\,\pi+\tfrac{5\pi}{4\sqrt{3}}+\sqrt{\tfrac{5}{7}}\,\beta\big),
\\
&\alpha_{13}=\big(\tfrac{\pi}{2\sqrt{3}}-\tfrac{2}{\sqrt{35}}\,\pi+\tfrac{2}{\sqrt{35}}\beta\big),
\quad
\alpha_{23}=\tfrac{1}{20}\big(5\sqrt{3}\,\pi+2\sqrt{35}\,\pi-2\sqrt{35}\,\beta\big).
\end{split}\]
For simplicity the constants $c_i{_\pm}$ and $\hat{\cU}$ were
equated to zero. One can check that after replacement of  $\log
(0_{+})$ by $\log(x)$ and substituting the final result into
exponenial formula $(23)$ in \cite{Acus2022} and then computing
the limit $x\to 0$ we recover the initial MV. To make the
verification simple when  $c_i{_\pm}$ and $\hat{\cU}$ are
included, one may choose concrete values for arbitrary free
constants $c_i{_\pm}$ and arbitrary unit bivector
$\hat{\cU}^2=-1$.

\section{Roots and arbitrary powers of MV}
\label{sec:ArbitraryPowers}

If GA logarithm is known, the powers of a MV  may be computed with
$\m{A}^r=\exp\bigl(r\log{\m{A}}\bigr)$, i.e.,  by multiplying
logarithm by  power value $r$,which may be either an integer or a
rational number, and then computing the exponential. In the
preprint~\cite{AcusDargysPreprint2020} we provided the algorithm
how to obtain all possible square roots ($r=1/2$) of MV for all
$n=3$ Clifford algebras. Here we want to show that the roots
presented in \cite{AcusDargysPreprint2020} as numerical examples
of algorithm  are consistent with the above exp-log formula, thus
actually we perform a cross check of 3D GA logarithm formulas by
different methods. It should be stressed that the logarithm
formula allows to find only a single\footnote{More precisely two
(plus/minus) roots.} square root, although there may exist, as
shown in ~\cite{AcusDargysPreprint2020}, many (up to 16 in case of
\cl{2}{1} algebra) roots. Thus, the GA logarithm function is not
universal enough, although it may be sometimes useful if only a
single fractional root, $r=1/n$ and $n\in\bbN$, is needed.

\vspace{3mm}
 \textbf{Example 12.} \cl{3}{0}. {\it Example 1 from
\cite{AcusDargysPreprint2020}.} Theorem~\ref{allCasesCl30} is used
to calculate  the root of  MV $\m{A}=\e{1}-2\e{12}$. In
$\cl{3}{0}$ algebra the logarithm is $\log \m{A}=\frac{\log
5}{2}-\frac{1}{2}\pi\e{23}+\arctan(\tfrac{1}{2})I$. Then the
square root $\sqrt{\m{A}}$ is
\begin{align*}
  \exp\bigl(\tfrac{1}{2}\log\m{A}\bigr)=&\frac{\sqrt[4]{5}}{\sqrt{2}}\bigl(\cos(\tfrac{1}{2}\arctan\tfrac{1}{2})\bigl(1-\e{23}\bigr)+\sin(\tfrac{1}{2}\arctan\tfrac{1}{2})\bigl(\e{1}+I\bigr) \bigr),
\end{align*}
which after simplification coincides with root $A_3$ in Example~1
in~\cite{AcusDargysPreprint2020}.

\vspace{4mm} \textbf{Example 13.}  \cl{3}{0}.  {\it Example 2 from
\cite{AcusDargysPreprint2020}.} Logarithm of MV
$\m{A}=-1+\e{3}-\e{12}+\tfrac{1}{2}I$ in $\cl{3}{0}$ is
\begin{align*}
  \log \m{A}=&\log\bigl(\tfrac{\sqrt{5}}{2}\bigr)-\tfrac{\log 5}{2}\e{3}+\tfrac{1}{2}\bigl(\pi-\arctan\tfrac{4}{3}\bigr)\e{12}
  +\bigl(-\pi+\arctan\tfrac{1}{2}\bigr)I ,
\end{align*}
Multiplication by $\frac{1}{2}$ and exponentiation gives the root
$A_3=\sqrt{\m{A}}=\frac{1}{2}\bigl(\e{3}+\e{12}\bigr)-I$ which
coincides with~\cite{AcusDargysPreprint2020}.

\textbf{Example 14.}  \cl{3}{0}. {\it Example 3 from
\cite{AcusDargysPreprint2020}.} Similarly, the logarithm of MV
$\m{A}=-1+\e{123}$  in $\cl{3}{0}$ is found to be $\log
\m{A}=\frac{\log 2}{2}+\frac{3}{4}\pi I$. Multiplication by
$\frac{1}{2}$ and exponentiation gives the root $A_1$ of Example~3
\cite{AcusDargysPreprint2020},
$\sqrt{\m{A}}=2^{1/4}\big(\cos\tfrac{3\pi}{8}+I\sin\tfrac{3\pi}{8}\big)=
\sqrt{-\tfrac{1}{2}+\tfrac{1}{\sqrt{2}}}+I\sqrt{\tfrac{1}{2}+\tfrac{1}{\sqrt{2}}}$.
 Likewise, an attempt to
compute the logarithm of $\m{A}=\e{1}+\e{12}$ yields empty set,
i.e., the logarithm and as a result the square root do not exist.

\textbf{Example 15.} \cl{3}{0}.  {\it Quaternion. Example 4 from
\cite{AcusDargysPreprint2020}.} In $\cl{3}{0}$ algebra the even MV
$\m{A}=1+\e{12}-\e{13}+\e{23}$ is equivalent to Hamilton
quaternion. The logarithm is $\log \m{A}=\log
2+\frac{\pi}{3\sqrt{3}}\bigl(\e{12}-\e{13}+\e{23}\bigr)$.
Multiplication by $\frac{1}{2}$ and exponentiation give the root
$A_3$ in Example~4 \cite{AcusDargysPreprint2020},
$\sqrt{\m{A}}=\tfrac{1}{\sqrt{6}}(3+\e{12}-\e{13}+\e{23})$.

\textbf{Example 16.}  \cl{0}{3}. {\it Example 6 from
\cite{AcusDargysPreprint2020}.} To compute the logarithm of MV
$\m{A}=\e{1}-2\e{23}$ in $\cl{0}{3}$ the
Theorem~\ref{allCasesCl03} was applied which gives  $\log
\m{A}=\frac{\log 3}{2} -\frac{\pi}{2}\e{23}+\frac{\log 3}{2} I$.
Multiplication by $\frac{1}{2}$ and exponentiation gives the root
$A_3$ of Example~6 in \cite{AcusDargysPreprint2020},
$\sqrt{\m{A}}=\tfrac12(d_2+d_1\e{1}-d_2\e{23}+I/d_2)$, where
$d_1=\sqrt{2-\sqrt{3}}$ and $d_2=\sqrt{2+\sqrt{3}}$.

\textbf{Example 17.} \cl{0}{3}.  {\it Example 7 from
\cite{AcusDargysPreprint2020}.} The logarithm of MV
$\m{A}=-\e{3}+\e{12}+4 I$ in $\cl{0}{3}$ algebra  is computed by
Theorem~\ref{allCasesCl03}. The result is $\log \m{A}=\frac{\log
320}{4} -\frac{1}{2}\arctan(\frac{1}{2})\e{3}
+\frac{1}{2}\arctan(\frac{1}{2})\e{12}
+\frac{\pi}{2}\bigl(1-I\bigr)\hat{\bu}+\frac{1}{4}\log\frac{5}{4}
I$. We have assumed that the discrete free constants are equal to
zero and retained only free unit vector $\hat{\bu}$ that satisfies
$\hat{\bu}^2=-1$. Multiplication by $\frac{1}{2}$ and
exponentiation then gives $\sqrt{\m{A}}=\frac{1}{2}c_2
+\frac{1}{2} c_1\bigl(\e{12}-\e{3}\bigr) + \frac{1}{2}c_2 I
+\bigl(1-I\bigr)\hat{\bu}$, where $c_1=\sqrt{\sqrt{5}-2}$ and
$c_2=\sqrt{2+\sqrt{5}}$, and corresponds to $A_3$ root in
Example~7 in \cite{AcusDargysPreprint2020}.  In particular, in
order to obtain the numerical value corresponding to
$V_2=\frac{1}{2}, V_3=0$ of \cite{AcusDargysPreprint2020} we have
to take $\hat{\bu}=-\frac{1}{2} \sqrt{5-\sqrt{5}-2
c_1}\e{1}+\frac{1}{2} (-1-c_1)\e{3}$.

\textbf{Example 18.} \cl{2}{1}.  {\it Example 8 from
\cite{AcusDargysPreprint2020}.} With the Theorem~\ref{theorem3}
one may ascertain that the logarithm of MV $\m{A}=\e{1}-2\e{23}$
in $\cl{2}{1}$ algebra does not exist what is in agreement with
square root absence of this MV in \cl{2}{1}. On the other hand,
the logarithm of MV $\m{A}=2+\e{1}+\e{13}$  is $\log
\m{A}=\frac{\log 2}{2}
-\frac{1}{\sqrt{2}}\arctanh\bigl(\frac{1}{\sqrt{2}}\bigr)
\bigl(\e{1}+\e{13}\bigr)$. Multiplication by $\frac{1}{2}$ and
exponentiation gives the root $A_5$ of example~8 in
\cite{AcusDargysPreprint2020}, $\sqrt{\m{A}}=\frac{1}{2}
\left(\sqrt{2-\sqrt{2}}\,(\e{1}+\e{13})+\sqrt{2
\left(2+\sqrt{2}\right)}\right)$.

Of course, after multiplication of logarithm by any integer or
rational number and subsequent exponentiation we can obtain  a
corresponding power of the MV. For example, in $\cl{3}{0}$ the
logarithm of  $\m{A}=\e{1}$ is $\log \m{A}=\frac{\pi}{2}\e{1}$.
Then, it is easy to check that after multiplication by
$\frac{1}{3}$ and exponentiation we obtain the cubic root
$\sqrt[3]{\e{1}}=\frac{1}{2} \bigl(\sqrt{3}+\e{1}\bigr)$.

\section{Relations of the logarithm to GA inverse trigonometric and hyperbolic functions}
\label{sec:relations}

Just like trigonometric and hyperbolic functions can be expressed
by exponentials (Euler and de Moivre formulas), the inverse
hyperbolic functions may be defined in terms of logarithms.
Therefore, in GA we can use the following definitions to compute
inverse hyperbolic and trigonometric functions of MV argument
$\m{A}$.\\ For hyperbolic inverse functions:
\begin{align}
\arctanh \m{A}=&\tfrac{1}{2}\bigl(\log(1 + \m{A}) - \log(1 - \m{A})\bigr),\\
  \arccoth \m{A}=&\begin{cases}
    \frac{1}{2}\bigl(\log(1 + \m{A}^{-1}) - \log(1 - \m{A}^{-1})\bigr),&\m{A}\neq 0,\\
    \frac{\pi}{2} I,& \m{A}= 0,
  \end{cases}\allowdisplaybreaks\\
\arccosh\m{A} =& \log\bigl(\m{A} + \sqrt{\m{A} - 1}\,\sqrt{\m{A} + 1}\bigr),\\
\arcsinh\m{A}= &\log\bigl(\m{A} + \sqrt{\m{A}^2 + 1}\bigr).
\end{align}
For inverse trigonometric functions:
\begin{align}
\arcsin\m{A}= &-I \log\bigl(\m{A} I + \sqrt{1-\m{A}^2}\bigr)\allowdisplaybreaks,\\
\arccos\m{A}= &\frac{\pi}{2}+I \log\bigl(\m{A} I + \sqrt{1-\m{A}^2}\bigr),\\
  \arctan \m{A}=&\frac{I}{2}\bigl(\log(1 - I \m{A}) - \log(1 + I \m{A})\bigr)\allowdisplaybreaks\label{invtrigFarctan},\\
\arccot \m{A}=&\begin{cases}
    \frac{1}{2}I\bigl(\log(1 -I \m{A}^{-1}) - \log(1 +I  \m{A}^{-1})\bigr),&\m{A}\neq 0,\\
    \frac{\pi}{2},& \m{A}= 0.
  \end{cases}\label{invtrigFarcCot}
\end{align}
These formulas are similar to those in the theory of real and
complex functions except that instead of the imaginary unit
the pseudoscalar appears in trigonometric functions. However,
earlier we have found~\cite{AcusDargysPreprint2020} that  in GA
the functions with the square root, in general, are multi-valued.
Thus at a first sight it may appear that the listed above
equations with square root are not valid in all circumstances.
Nonetheless, our preliminary numerical experiments show that they,
in fact, are satisfied for all possible individual
plus/minus pairs of square roots\footnote{This property does not
allow us to write the equality sign between GA general
expression $\log\sqrt{\m{B}}$ and $\frac12\log\m{B}$.} (see
Example~8 in~\cite{Acus2022} and Example~19 below in this
section).

With the above formulas for hyperbolic and trigonometric functions
one can construct the following identities for generic
MVs:\footnote{Trigonometric functions are defined only for
algebras where the pseudoscalar 1) belongs to a center of an algebra, i.e. commutes
commutative with remaining elements and 2) satisfy $I^2=-1$. In the
considered 3D algebras only for $\cl{3}{0}$ and $\cl{1}{2}$.}
\begin{align}
\sinh\m{A}= &\tfrac{1}{2} \bigl(\exp(\m{A}) - \exp(-\m{A})\bigr),\\
\cosh\m{A}= &\tfrac{1}{2} \bigl(\exp(\m{A}) + \exp(-\m{A})\bigr)\allowdisplaybreaks,\\
\tanh\m{A}=&\sinh\m{A}(\cosh\m{A})^{-1}=\bigl(\exp(\m{A}) - \exp(-\m{A})\bigr) \bigl(\exp(\m{A}) + \exp(-\m{A})\bigr)^{-1},\\
\coth\m{A}=&\cosh\m{A}(\sinh\m{A})^{-1}=\bigl(\exp(\m{A}) +
\exp(-\m{A})\bigr) \bigl(\exp(\m{A}) -
\exp(-\m{A})\bigr)^{-1}\allowdisplaybreaks.
\end{align}
\begin{align}
\sin\m{A}= &\tfrac{1}{2}I^{-1} \bigl(\exp(I \m{A}) - \exp(-I \m{A})\bigr),\\
\cos\m{A}= &\tfrac{1}{2} \bigl(\exp(I \m{A}) + \exp(-I \m{A})\bigr)\allowdisplaybreaks,\\
\tan\m{A}=&\sin\m{A}(\cos\m{A})^{-1}=-I\bigl(\exp(I\m{A}) - \exp(-I\m{A})\bigr) \bigl(\exp(I\m{A}) + \exp(-I\m{A})\bigr)^{-1},\\
\cot\m{A}=&\cos\m{A}(\sin\m{A})^{-1}=I \bigl(\exp(I \m{A}) +
\exp(-I \m{A})\bigr) \bigl(\exp(I \m{A}) - \exp(-I
\m{A})\bigr)^{-1}.
\end{align}
We have not investigated how the presented formulas work in  case
when the MV square root or logarithm allows answer that depends on
non-discrete free parameters and when the inverse MVs can't be be
computed. Also, we have not considered MV logarithms that allow
infinite coefficients at some of basis MVs.

\vspace{3mm}
 \textbf{Example 19.} {\it Inverse MV hyperbolic
functions.} To save space we will restrict ourselves to numerical
examples only for $\cl{3}{0}$ generic MV $\m{A}=-1-5 \e{1}+7
\e{2}-9 \e{3}+7 \e{12}-5 \e{13}+9 \e{23}+9 I$. Then we find the
following inverse hyperbolic functions,
\begin{align*}
\arctanh \m{A}=&\kern1em\begin{aligned}[t]
    &0.0544776 &\kern-1em-&0.0683983 \e{1}&\kern-1em-&0.0034179 \e{2}&\kern-1em+&0.0712752 \e{3}\\
  \kern-1em-&0.0259578 \e{12}&\kern-1em-&0.0571283 \e{13}&\kern-1em+&0.0036554 \e{23}&\kern-1em+&1.5447402
  I,
  \end{aligned}\\
\arccoth \m{A}=&\kern1em\begin{aligned}[t]
 & 0.0544776&\kern-1em-&0.0683983 \e{1}&\kern-1em-&0.0034179 \e{2}&\kern-1em-&0.0712752 \e{3}\\
\kern-1em-&0.0259578 \e{12}&\kern-1em-&0.0571283
\e{13}&\kern-1em+&0.0036555 \e{23}&\kern-1em-&0.0260523 I,
  \end{aligned}\\
 \arccosh \m{A}=&\kern1em\begin{aligned}[t]
   &3.1995844 &\kern-1em+&0.6349751 \e{1}&\kern-1em+&0.6477695 \e{2}&\kern-1em+&0.3396621 \e{3}\\
  \kern-1em+&0.9970274 \e{12}&\kern-1em+&0.4603461 \e{13}&\kern-1em+&0.7081115 \e{23}&\kern-1em+&1.0647020
  I.
  \end{aligned}
\end{align*}
For identities that contain square roots $\sqrt{\m{A}\pm 1}$, for
example  $\arccosh \m{A}$ or $\arcsinh \m{A}$, all four
roots are valid. Below they have been calculated by algorithm
described in~\cite{AcusDargysPreprint2020},
\begin{align*}
\textit{Root 1 and 2}:\\
\sqrt{\m{A} - 1}=\pm(&\kern1em\begin{aligned}[t]
\kern-1em-&2.3936546&\kern-1em-&0.3144420 \e{1}&\kern-1em-&1.3708134 \e{2}&\kern-1em+&0.3806804\e{3}\\
\kern-1em-&1.7824116 \e{12}&\kern-1em-&0.1086429 \e{13}&\kern-1em-&1.6154750 \e{23}&\kern-1em-&2.0134421 I),
 \end{aligned}\\
  \indent \textit{Root 3 and 4}:\\
\sqrt{\m{A} - 1}=\pm(&\kern1em\begin{aligned}[t]
\kern-1em-&0.1660207 &\kern-1em+&   2.4324037 \e{1}&\kern-1em+&1.1337774 \e{2}&\kern-1em+&2.0055931 \e{3}\\
\kern-1em+&2.1654007 \e{12}&\kern-1em+&1.9165921
\e{13}&\kern-1em+&1.0892769 \e{23}&\kern-1em+&1.9243691 I).
  \end{aligned}
\end{align*}
And similarly for
\begin{align*}
\sqrt{\m{A} + 1}=\pm\{&\kern1em\begin{aligned}[t]
\kern-1em-&2.6330243&\kern-1em-&0.1908183 \e{1}&\kern-1em-&1.3218252 \e{2}&\kern-1em+&0.4871255 \e{3}\\
\kern-1em-&1.6829550 \e{12}&\kern-1em-&0.0102534 \e{13}&\kern-1em-&1.5705147 \e{23}&\kern-1em-&1.9117486 I\},\\[5pt]
 \end{aligned}\\
\sqrt{\m{A} + 1}=\pm\{&\kern1em\begin{aligned}[t]
\kern-1em-&0.2910283 &\kern-1em+& 2.6047118 \e{1}&\kern-1em+&1.0343473 \e{2}&\kern-1em+&2.2416242 \e{3}\\
\kern-1em+&2.0981982 \e{12}&\kern-1em+&2.0727864
\e{13}&\kern-1em+&0.9499284 \e{23}&\kern-1em+&1.8337753 I\}.
  \end{aligned}
\end{align*}
It is important to stress that, in general, the individual
formulas ($\arccos \m{A}$, $\arcsin \m{A}$ and their hyperbolic
analogues) that contain the sets of  roots  yield different
function values for four different roots in the above
listed sets.
\begin{align*}
\arcsinh \m{A}=\pm(&\kern1em\begin{aligned}[t]
\kern-1em\ &3.2035891&\kern-1em+&0.6313828 \e{1}&\kern-1em+&0.6490577 \e{2}&\kern-1em+&0.3351515 \e{3}\\
\kern-1em+&0.9974654 \e{12}&\kern-1em+&0.4571790 \e{13}&\kern-1em+&0.7100715 \e{23}&\kern-1em+&1.0647010 I).\\
\end{aligned}
\end{align*}
\begin{align*}
\arcsinh \m{A}=\pm(&\kern1em\begin{aligned}[t]
\kern-1em\ &0.4835482 &\kern-1em+& 2.5061943 \e{1}&\kern-1em-&0.7303556 \e{2}&\kern-1em+&3.0588480 \e{3}\\
\kern-1em-&0.0989201 \e{12}&\kern-1em+&2.1904765
\e{13}&\kern-1em-&1.1645414 \e{23}&\kern-1em+&2.5756463 I).
  \end{aligned}
\end{align*}

\textbf{Example 20.} {\it Inverse trigonometric functions of MV.}
Numerical answers for $\cl{3}{0}$ generic MV $\m{A}=-1-5 \e{1}+7
\e{2}-9 \e{3}+7 \e{12}-5 \e{13}+9 \e{23}+9 I$ in a form of list
for roots 1-4,
\begin{align*}
\arcsin \m{A}=\{&\kern1em\begin{aligned}[t]
\kern-1em\ &2.5745928&\kern-1em+&0.1233316 \e{1}&\kern-1em-&2.3715122 \e{2}&\kern-1em+&1.3713947 \e{3}\\
\kern-1em-&2.8712504 \e{12}&\kern-1em+&0.3732007 \e{13}&\kern-1em-&2.8706092 \e{23}&\kern-1em-&0.4882339 I,\\
\kern-1em\ &2.6354984 &\kern-1em+& 0.7081116 \e{1}&\kern-1em-&0.4603462 \e{2}&\kern-1em+&0.9970274 \e{3}\\
\kern-1em-&0.3396621 \e{12}&\kern-1em+&0.6477695 \e{13}&\kern-1em-&0.6349751 \e{23}&\kern-1em-&3.1995845 I,\\
\kern-1em+&0.5669998  &\kern-1em-&0.1233316 \e{1}&\kern-1em+&2.3715122 \e{2}&\kern-1em-&1.3713947 \e{3}\\
\kern-1em+&2.8712504 \e{12}&\kern-1em-&0.3732007 \e{13}&\kern-1em+&2.8706092 \e{23}&\kern-1em+&0.4882339 I,\\
\kern-1em\ &0.5060943  &\kern-1em-&0.7081116 \e{1}&\kern-1em+&0.4603462 \e{2}&\kern-1em-&0.9970274 \e{3}\\
\kern-1em+&0.3396621 \e{12}&\kern-1em-&0.6477695
\e{13}&\kern-1em+&0.6349751 \e{23}&\kern-1em+&3.1995845 I\}.
  \end{aligned}
\end{align*}
Since formulas for arc sine and cosine also include square roots
we obtain four different values for these functions too,
\begin{align*}
\arccos \m{A}=\pm\{&\kern1em\begin{aligned}[t]
\kern-1em-&1.0037965&\kern-1em-&0.1233316 \e{1}&\kern-1em+&2.3715122 \e{2}&\kern-1em-&1.3713947 \e{3}\\
\kern-1em+&2.8712504 \e{12}&\kern-1em-&0.3732007 \e{13}&\kern-1em+&2.8706092 \e{23}&\kern-1em+&0.4882339 I,\\
\kern-1em-&1.0647021  &\kern-1em-&0.7081116 \e{1}&\kern-1em+&0.4603462 \e{2}&\kern-1em-&0.9970274 \e{3}\\
\kern-1em+&0.3396621 \e{12}&\kern-1em-&0.6477695
\e{13}&\kern-1em+&0.6349751 \e{23}&\kern-1em+&3.1995845 I\}.
\end{aligned}
\end{align*}
On the other hand the  trigonometric tangent and cotangent have a
single value  since the square root here is absent,
Eqs~\eqref{invtrigFarctan} and \eqref{invtrigFarcCot},
\begin{align*}
 \arctan \m{A}=&\kern1em\begin{aligned}[t]
\kern-1em\ &1.5171201 &\kern-1em+&  0.0678435 \e{1}&\kern-1em+&0.0036019 \e{2}&\kern-1em+&0.0705863 \e{3}\\
\kern-1em+&0.0260071 \e{12}&\kern-1em+&0.0566409
\e{13}&\kern-1em-&0.0033708 \e{23}&\kern-1em+&0.0259164 I,
\end{aligned}
\end{align*}
\begin{align*}
 \arccot \m{A}=&\kern1em\begin{aligned}[t]
\kern-1em\ &0.0536762 &\kern-1em-&0.0678435 \e{1}&\kern-1em-&0.0036019 \e{2}&\kern-1em-&0.0705863 \e{3}\\
\kern-1em-&0.0260071 \e{12}&\kern-1em-&0.0566409
\e{13}&\kern-1em+&0.0033708 \e{23}&\kern-1em-&0.0259164 I.
 \end{aligned}
\end{align*}

\section{Discussion and conclusions}
\label{sec:discussion}

The  logarithm together with the
exponential~\cite{Acus2022,Dargys2022a} and square
root~\cite{AcusDargysPreprint2020}  are the most important
functions in Clifford geometric algebra (GA). Starting from the
respective exponential functions we presented  here, as far
as we know,  for the first time the basis-free formulas for
logarithms in all 3D GAs. The  formulas for both the generic and
special  cases may be directly applied in GA programming. They
were cross-checked using the basis-free GA exponential functions
found in~\cite{Acus2022}. The derived formulas were implemented in
\textit{Mathematica} and tested with thousands of randomly
generated multivectors~\cite{AcusDargys2017}.  In all cases the
exponentiation of the logarithm was found to simplify to the
initial MV.

Using numerical experiments~\cite{AcusDargys2017} we observed
that, in accord with the suggestion in \cite{Chappell2015}, the
principal value of the logarithm can be defined as a GA logarithm
having the smallest determinant norm. In almost all cases the
principal MV logarithm  is attained by setting arbitrary integer
parameters $c_i$ in generic logarithms
(Theorems~\ref{log03coordfree}, \ref{allCasesCl30},
\ref{log21thmTheorem}) to zero. Exceptions from this rule,
however, may occur in the case of simple specific MVs, for which
commuting MVs may exist ( Secs.~\ref{freeParameters} and
\ref{freeParametersCl30}), and therefore not  restricted  by free
MVs (Eqs~\eqref{periodicityCL03} and \eqref{periodicityCL30CL12}).
Apart from discrete parameters $c_i$, we have also found that
continuous parameters represented by free unit vectors $\hat{\bu}$
or bivectors $\hat{\cU}$ may be included in special cases as well.
The parameters vanish after exponentiation of the logarithm and do
not contribute to the MV norm. Recently we have found that
such free parameters may be also introduced into lower
dimensional, quaternionic-type Clifford
algebras~\cite{Dargys2022Log}. However, more investigations are
needed in this direction.

Also,  the relation between the GA logarithm and square root of MV
was investigated. The known formula $\sqrt{\m{A}}=\exp
\bigl(\frac12\log (\m{A})\bigr)$ served as an additional check
correctness of GA logarithms. Unfortunately, the formula allows to
compute only a single square root from many possible roots that
may exist in GA~\cite{AcusDargysPreprint2020}. Nevertheless, such
a comparison was found to be very useful for testing purposes.
Indeed, a test of square root of a MV is an algebraic problem
since it reduces to a solution of system of algebraic equations.
On the other hand, inversion of exponential used in finding the GA
logarithm in the present paper requires solving a system of
transcendental equations (Appendix~\ref{coordFormCL03}), a problem
which is much more difficult (but at the same time more general)
task. The mentioned exp-log relation also allows to check the
condition whether the MV logarithm exists at all. Indeed, since we
know how to calculate GA exponential~\cite{AcusDargysPreprint2020}
of arbitrary MV multiplied by factor~$\frac12$, from this follows
that it is $\log (\m{A})$ function which determines the existence
condition for $\sqrt{\m{A}}$ to exist. As a test, we have checked
using our algorithm~\cite{AcusDargysPreprint2020} that for each MV
indeed there exists a single square root that is in  agreement
with the identity $\exp \bigl(\frac12\log
(\m{A})\bigl)=\sqrt{\m{A}}$.

In conclusion, in the present paper the basis-free expressions
have been found for GA logarithms in all 3D real algebras. The
logarithm was found to exist for all MVs in case of real \cl{0}{3}
algebra. In Clifford algebra \cl{3}{0} (and \cl{1}{2}) the
logarithm exists for almost all MVs, except very small MV class
which satisfies the condition $(a_{+}^2 + a_{-}^2 = 0)\land
(a_{0}^2 + a_{123}^2= 0)$.  For example, the logarithm of MV
$\e{1}\pm\e{12}$ cannot be computed in Euclidean \cl{3}{0}
algebra. On the other hand in \cl{2}{1} algebra the GA logarithm
is absent in large sectors of a real coefficient space.


%
\bibliographystyle{REPORT}
\bibliography{logarithm3D}

\appendix

\section{\label{coordFormCL03}Outline of derivation of generic logarithm formula
in a coordinate form in \cl{0}{3}}

The MV exponential  $\exp\m{A}=\m{B}=b_0+\bb+\cB+b_{123}I$, where
$\bb=b_1\e{1}+b_2\e{2}+b_3\e{3}$ and
$\cB=b_{12}\e{12}+b_{13}\e{13}+b_{23}\e{23}$, in the coordinate
form was constructed in~\cite{Dargys2022a}.  For completeness,
below we reproduce the expressions for scalar coefficients of
$\m{B}$,
\begin{align}
&\begin{aligned}
\phantom{{}_{99}}b_{0}&=\tfrac{1}{2}\ee^{a_{0}} \bigl(\ee^{a_{123}}\cos a_{+} + \ee^{-a_{123}} \cos a_{-}\bigr),\\
b_{123}&=\tfrac{1}{2}\ee^{a_{0}} \bigl(\ee^{a_{123}}\cos a_{+}
-\ee^{-a_{123}}\cos a_{-}\bigr),
\end{aligned}\label{exp03components07}
\allowdisplaybreaks \\[4pt]
&\begin{aligned}
  \phantom{{}_{99}}b_{1}&=\tfrac{1}{2}\ee^{a_{0}} \Bigl(\ee^{a_{123}}(a_{1}-a_{23})\frac{\sin a_{+}}{a_{+}} + \ee^{-a_{123}}(a_{1}+a_{23})\frac{\sin a_{-}}{a_{-}}\Bigr),\\
\phantom{{}_{99}}b_{2}&=\tfrac{1}{2}\ee^{a_{0}} \Bigl(\ee^{a_{123}}(a_{2}+a_{13})\frac{\sin a_{+}}{a_{+}} +\ee^{-a_{123}}(a_{2}-a_{13})\frac{\sin a_{-}}{a_{-}}\Bigr),\\
\phantom{{}_{99}}b_{3}&=\tfrac{1}{2}\ee^{a_{0}} \Bigl(\ee^{a_{123}}(a_{3}-a_{12})\frac{\sin a_{+}}{a_{+}} +\ee^{-a_{123}}(a_{3}+a_{12})\frac{\sin a_{-}}{a_{-}}\Bigr),
\end{aligned}\label{exp03componentsvec}\allowdisplaybreaks\\[4pt]
&
\begin{aligned}
\phantom{{}_{9}}b_{12}&=\tfrac{1}{2}\ee^{a_{0}} \Bigl(-\ee^{a_{123}}(a_{3}-a_{12})\frac{\sin a_{+}}{a_{+}} +\ee^{-a_{123}}(a_{3}+a_{12})\frac{\sin a_{-}}{a_{-}}\Bigr),\\
\phantom{{}_{9}}b_{13}&=\tfrac{1}{2}\ee^{a_{0}} \Bigl(\ee^{a_{123}}(a_{2}+a_{13})\frac{\sin a_{+}}{a_{+}} -\ee^{-a_{123}}(a_{2}-a_{13})\frac{\sin a_{-}}{a_{-}}\Bigr),\\
\phantom{{}_{9}}b_{23}&=\tfrac{1}{2}\ee^{a_{0}} \Bigl(-\ee^{a_{123}}(a_{1}-a_{23})\frac{\sin a_{+}}{a_{+}} +\ee^{-a_{123}}(a_{1}+a_{23})\frac{\sin a_{-}}{a_{-}}\Bigr),
\end{aligned}\label{exp03componentsbiv}
\end{align}
where  $a_{-}$ and $a_{+}$ are defined in Eqs~\eqref{aminusCl03}
and~\eqref{apliusCl03}, respectively. Explicit generic formula for
a MV logarithm $\text{log}\m{B}=\m{A}$ can be derived by inverting
the equations
\eqref{exp03components07}-\eqref{exp03componentsbiv}. To this end,
using \eqref{exp03components07}-\eqref{exp03componentsbiv} we will
construct a system of eight nonlinear trigonometric equations for
coefficients at corresponding basis elements. Since two of the
coefficients, $a_0$ and $a_{123}$,  stand alone, i.e., they are
not associated with remaining ones in \eqref{aminusCl03}
and~\eqref{apliusCl03}, as a first step in this procedure we
partially solve a pair of equations at scalar and pseudoscalar, $
b_{0}$ and $b_{123}$, to get two new functions
\begin{subequations}
\begin{align}
  a_{0}=&g^\prime_{0}(b_{0},b_{123},a_{-},a_{+}),\label{loggenerala0}\\
  a_{123}=&g^\prime_{123}(b_{0},b_{123},a_{-},a_{+}),\label{loggenerala7}
\end{align}
\end{subequations}
where the pair $(a_{-},a_{+})$ includes only vector and bivector
components as \eqref{aminusCl03} and Eqs~\eqref{apliusCl03} show.
After substituting new $a_{0}$ and $a_{123}$ into remaining
equations \eqref{exp03componentsvec} and
\eqref{exp03componentsbiv} we obtain a system of nonlinear
trigonometric equations for coefficients at vector and bivector
only which, unfortunately, cannot be solved by computer algebra
system  in this form. Further progress can be achieved by
observing that the following relation between the coefficients
holds in the generic case
\begin{subequations}
\begin{align}
  a_{-}=&\arctan(b_{0}-b_{123},b_{-}),\label{logaPlusbPlusrelationsam}\\
  a_{+}=&\arctan(b_{0}+b_{123},b_{+}),\label{logaPlusbPlusrelationsap}
\end{align}
\end{subequations}
where the pair $(b_{-},b_{+})$ is given by Eqs \eqref{apliusCl03}
and \eqref{aminusCl03} after replacement of $\ba$ and $\cA$, by
vector $\bb$ and bivector $\cB$, respectively. Then, replacing all
occurrences of $(a_{-},a_{+})$ by \eqref{logaPlusbPlusrelationsam}
and \eqref{logaPlusbPlusrelationsap} in \eqref{exp03componentsvec}
and \eqref{exp03componentsbiv} (with already replaced $a_{0}$ and
$a_{123}$ as was described above) we obtain a linear system for
coefficients $a_{i},a_{ij}$, that can be solved in a
straightforward way. Finally, after substituting $a_{-}$ and
$a_{+}$ as given by \eqref{logaPlusbPlusrelationsam} and
\eqref{logaPlusbPlusrelationsap} into \eqref{loggenerala0}
and\eqref{loggenerala7} and performing simplifications we find the
scalar coefficients of $\log\m{B}=\m{A}=a_0+\ba+\cA+a_{123}I$
expressed in terms of two-argument arc tangent functions (see
Subsec.~\ref{expLogTan}),
\begin{align}\label{log03coord1}
a_{0}=& \tfrac{1}{4} \left(\log \left((b_{0}-b_{123})^2+b_ {-}^{2}\right)+\log \left((b_{0}+b_{123})^2+b_ {+}^{2}\right)\right),\notag\\
a_{123}=& \tfrac{1}{4} \left(\log \left((b_{0}+b_{123})^2+b_
{+}^{2}\right)-\log \left((b_{0}-b_{123})^2+b_
{-}^{2}\right)\right),\notag\\
\end{align}
\begin{align}\label{log03coord2}
  a_{1}=& \tfrac{1}{2} \left(\frac{(b_{1}+b_{23}) \arctan \left(b_{0}-b_{123},b_{-}\right)}{b_{-}}+\frac{(b_{1}-b_{23}) \arctan \left(b_{0}+b_{123},b_{+}\right)}{b_{+}}\right),\notag\\
  a_{2}=& \tfrac{1}{2} \left(\frac{(b_{2}-b_{13}) \arctan \left(b_{0}-b_{123},b_{-}\right)}{b_{-}}+\frac{(b_{2}+b_{13}) \arctan \left(b_{0}+b_{123},b_{+}\right)}{b_{+}}\right),\notag\\
  a_{3}= &\tfrac{1}{2} \left(\frac{(b_{3}+b_{12}) \arctan \left(b_{0}-b_{123},b_{-}\right)}{b_{-}}+\frac{(b_{3}-b_{12}) \arctan
  \left(b_{0}+b_{123},b_{+}\right)}{b_{+}}\right),\notag\\
\end{align}
\begin{align}\label{log03coord3}
a_{12}=&\tfrac{1}{2} \left(\frac{(b_{3}+b_{12}) \arctan \left(b_{0}-b_{123},b_{-}\right)}{b_{-}}+\frac{(b_{12}-b_{3}) \arctan \left(b_{0}+b_{123},b_{+}\right)}{b_{+}}\right),\notag\\
  a_{13}=& \tfrac{1}{2} \left(\frac{(b_{13}-b_{2}) \arctan \left(b_{0}-b_{123},b_{-}\right)}{b_{-}}+\frac{(b_{2}+b_{13}) \arctan \left(b_{0}+b_{123},b_{+}\right)}{b_{+}}\right),\notag\\
  a_{23}=&\tfrac{1}{2} \left(\frac{(b_{1}+b_{23}) \arctan \left(b_{0}-b_{123},b_{-}\right)}{b_{-}}+\frac{(b_{23}-b_{1}) \arctan
  \left(b_{0}+b_{123},b_{+}\right)}{b_{+}}\right).
\end{align}
Here $b_{+}$ and $b_{-}$ are given by \eqref{apliusCl03} and
\eqref{aminusCl03} after replacement of $\ba$ and $\cA$,
respectively, by vector $\bb$ and bivector $\cB$. The equations
\eqref{log03coord1}-\eqref{log03coord3} provide a generic solution
of the inverse problem for MV logarithm in a coordinate form. It
must be remembered that in the two-argument arc tangent functions
the argument order and properties follow {\it Mathematica}
convention (see Subsec.~\ref{expLogTan}). To have a general GA
logarithm formula of Theorem~\ref{theoremCl30} in a basis-free
form, it is enough to multiply the coefficients
in~\eqref{log03coord1}-\eqref{log03coord3} by respective basis
elements, to add and regroup the resulting GA expression into
coordinate-free form.

\end{document}